\documentclass[12pt]{article}
\usepackage{latexsym}
\usepackage{amsfonts}
\usepackage{amssymb}
\usepackage{amscd}
\usepackage{array}
\usepackage{amsmath}
\usepackage{epsfig}
\usepackage{a4}
\usepackage{longtable}

\begin{document}
\title{\textbf{Stringy $E$-functions of varieties with $A$-$D$-$E$ singularities}\\ \date{}}
\author{Jan Schepers\footnote{Research Assistant of the Fund for Scientific Research
- Flanders (Belgium) (F.W.O.), \textsc{Katholieke Universiteit
Leuven, Departement Wiskunde, Celestijnenlaan 200B, B-3001 Leuven,
Belgium}, \emph{E-mail address}:
jan.schepers@wis.kuleuven.ac.be.}}

\maketitle
\begin{center}
\footnotesize{\textbf{Abstract}}
\end{center}
{\footnotesize The stringy $E$-function for normal irreducible
complex varieties with at worst log terminal singularities was
introduced by Batyrev. It is defined by data from a log
resolution. If the variety is projective and Gorenstein and the stringy
$E$-function is a polynomial, Batyrev also defined the stringy
Hodge numbers as a generalization of the Hodge numbers of
nonsingular projective varieties, and conjectured that they are
nonnegative. We compute explicit formulae for the contribution of
an $A$-$D$-$E$ singularity to the stringy $E$-function in
arbitrary dimension. With these results we can say when the
stringy $E$-function of a variety with such singularities is a
polynomial and in that case we prove that the stringy Hodge
numbers are nonnegative.}
 \\
\section{Introduction}

\textbf{1.1.} In \cite{Batyrev}, Batyrev defined the stringy
$E$-function for normal irreducible complex algebraic varieties,
with at worst log terminal singularities. With this function he
was able to formulate a topological mirror symmetry test for
Calabi-Yau varieties with singularities. Before stating the
definition of the stringy $E$-function, we recall some other
definitions.

Let $X$ be a complex algebraic variety. One defines the
Hodge-Deligne polynomial $H(X;u,v)\in \mathbb{Z}[u,v]$ by
\[H(X;u,v) = \sum_{i=0}^{2d} (-1)^i \sum_{p,q}
h^{p,q}(H_c^i(X,\mathbb{C}))u^pv^q,\] where $h^{p,q}$ denotes the
dimension of the $(p,q)$-component of the mixed Hodge structure on
$H^i_c(X,\mathbb{C})$. A nice introduction to Deligne's mixed
Hodge theory and to this definition can be found in
\cite{Srinivas} (pay attention to the extra factor $(-1)^{p+q}$
that the author has inserted there). The Hodge-Deligne polynomial
is a generalized Euler characteristic, that is, it satisfies:
\begin{itemize}
\item[(1)]$H(X)=H(X\setminus Y)+H(Y)$ where $Y$ is Zariski-closed
in $X$,
\item[(2)]$H(X\times X')=H(X)\cdot H(X')$.
\end{itemize}
Note that $H(X;1,1)=\chi(X)$, the topological Euler characteristic
of $X$.\\

\noindent \textbf{1.2.} A normal irreducible complex variety $X$
is called $\mathbb{Q}$-Gorenstein if $rK_X$ is Cartier for some
$r\in\mathbb{Z}_{>0}$. Take a log resolution $\varphi:
\widetilde{X}\to X$ (i.e. a proper birational morphism from a
nonsingular variety $\widetilde{X}$ such that the exceptional
locus of $\varphi$ is a divisor whose components $D_1,\ldots,D_s$
are smooth and have normal crossings). Then we have
$rK_{\widetilde{X}}-\varphi^*(rK_X)=\sum_i b_i D_i$, with $b_i\in
\mathbb{Z}$. This is also formally written as
$K_{\widetilde{X}}-\varphi^*(K_X)=\sum_i a_iD_i$, where
$a_i=\frac{b_i}{r}$. The variety $X$ is called terminal,
canonical, log terminal and log canonical if $a_i>0,a_i \geq 0,a_i
> -1,a_i\geq -1,$ respectively, for all $i$ (this is independent of
the chosen log resolution). The difference
$K_{\widetilde{X}}-\varphi^*(K_X)$ is called the
\emph{discrepancy}.\\

\noindent \textbf{1.3.} Now we are ready to define the stringy
$E$-function. We discuss its properties and give the additional
definitions of the stringy Euler number and the stringy Hodge
numbers. All of this goes back to Batyrev \cite{Batyrev}.\\

\noindent \textbf{Definition.} Let $X$ be a normal irreducible
complex variety with at most log terminal singularities and let
$\varphi:\widetilde{X}\to X$ be a log resolution. Denote the
irreducible components of the exceptional locus by $D_i$, $i\in
I$, and write $D_J$ for $\cap_{j\in J} D_j$ and $D_J^{\circ}$ for
$D_J\setminus \cup_{j\in I\setminus J} D_j$, where $J$ is any
subset of $I$ ($D_{\emptyset}$ is taken to be $\widetilde{X}$).
The stringy $E$-function of $X$ is
\[E_{st}(X;u,v):= \sum_{J\subseteq I} H(D_J^{\circ};u,v) \prod_{j\in J}
\frac{uv-1}{(uv)^{a_j+1}-1}, \] where $a_j$ is the discrepancy
coefficient of $D_j$ and where the product $\prod_{j\in J}$ is 1
if $J=\emptyset$.\\

Batyrev proved that this definition is independent of the chosen
log resolution. His proof uses motivic integration. An overview of
this theory is provided in \cite{Veys}.\\

\noindent \textbf{Remark.}
 
\begin{itemize}
\item[(1)] If $X$ is smooth, then $E_{st}(X)=H(X)$ and if $X$
admits a crepant resolution $\varphi:\widetilde{X} \to X$ (i.e.
such that the discrepancy is 0), then
$E_{st}(X)=H(\widetilde{X})$.
\item[(2)] If $X$ is Gorenstein (i.e. $K_X$ is Cartier), then all
$a_i\in \mathbb{Z}_{\geq 0}$ and $E_{st}(X)$ becomes a rational
function in $u$ and $v$. It is then an element of
$\mathbb{Z}[[u,v]]\cap \mathbb{Q}(u,v)$.
\item[(3)] The \emph{stringy Euler number} of $X$ is defined as
\[\lim_{u,v \to 1} E_{st}(X;u,v)=\sum_{J\subseteq I}
\chi(D_J^{\circ})\prod_{j\in J} \frac{1}{a_j+1}.\]
\end{itemize}

\noindent \textbf{1.4.} Assume moreover that $X$ is projective of
dimension $d$. Then Batyrev proved the following instance of
Poincar\'e and Serre duality:
\begin{itemize}
\item[(i)] $E_{st}(X;u,v)=(uv)^dE_{st}(X;u^{-1},v^{-1}),$
\item[(ii)]$E_{st}(X;0,0)=1$.
\end{itemize}
If $X$ has at worst Gorenstein canonical singularities and if
$E_{st}(X;u,v)$ is a polynomial $\sum_{p,q} a_{p,q}u^pv^q$, he
defined the \emph{stringy Hodge numbers} of $X$ as
$h_{st}^{p,q}(X):=(-1)^{p+q}a_{p,q}$. It is clear that
\begin{itemize}
\item[(1)] they can only be nonzero for $0\leq p,q\leq d$,
\item[(2)] $h_{st}^{0,0}=h_{st}^{d,d}=1$,
\item[(3)] $h_{st}^{p,q}=h_{st}^{q,p}=h_{st}^{d-p,d-q}=h_{st}^{d-q,d-p}$,
\item[(4)] if $X$ is smooth, the stringy Hodge numbers are equal
to the usual Hodge numbers.
\end{itemize}

\noindent \textbf{Conjecture (Batyrev).}  \textsl{The stringy Hodge
numbers are nonnegative.}\\

\noindent \textbf{Example.} The conjecture is true for varieties that admit a crepant
resolution. This is the case for all canonical surface
singularities, which are exactly the two-dimensional $A$-$D$-$E$
singularities \cite[p.375]{Reid} (see also Theorem 5.1 for $m=3$).\\

\noindent \textbf{Remark.} For a complete surface $X$ with at most log terminal
singularities, Veys showed that
\[E_{st}(X)=\sum_{p,q\in \mathbb{Z}}(-1)^{p+q} h_{st}^{p,q} u^pv^q
+\sum_{r\notin \mathbb{Z}} h_{st}^{r,r}(uv)^r, \] with all
$h_{st}^{p,q}$ and $h_{st}^{r,r}$ nonnegative \cite[p.138]{Veys2}.\\

\noindent \textbf{1.5.} In this paper, we will compute in
arbitrary dimension the contribution of an $A$-$D$-$E$ singularity
to the stringy $E$-function. This has already been done by Dais
and Roczen in the three-dimensional case (see \cite{DaisRoczen}),
but their computation of some discrepancy coefficients in the $D$
and $E$ cases is inaccurate and this leads to incorrect formulae
in these cases. We correct and considerably simplify their
formulae (also for type $A$). We construct a log resolution for
all higher dimensional $A$-$D$-$E$ singularities (based on the
calculation by Dais and Roczen of a log resolution for the
three-dimensional $A$-$D$-$E$'s), and again we are always able to
obtain a fairly simple formula for their stringy $E$-function. For
the contribution of an ($m-1$)-dimensional singularity of type
$D_n$ (where $m$ is odd and $n=2k$ is even) we find for example
\[ 1+\frac{(uv-1)}{((uv)^{(2k-1)(m-3)+1}-1)}\left(\sum_{i=1}^{2k-1}(uv)^{i(m-3)+1}
+(uv)^{k(m-3)+1}\right).   \] Then using our concrete formulae, we
can prove the following theorem.\\

\noindent \textbf{Theorem.}  \textsl{Let $X$ be a projective complex
variety of dimension at least 3 with at most $A$-$D$-$E$
singularities. The stringy $E$-function of $X$ is a polynomial if
and only if $X$ has dimension 3 and all singularities are of type
$A_n$ ($n$ odd) and/or $D_n$ ($n$ even). In that case, the stringy
Hodge numbers of $X$ are positive.}\\

In the next section we recall the definition of the $A$-$D$-$E$
singularities and we construct a log resolution for them. In
section 3 and 4, we compute the Hodge-Deligne polynomials and the
discrepancy coefficients that we need, respectively. In section 5
we give the resulting formulae and prove the theorem.\\

\section{$A$-$D$-$E$ singularities and their desingularization}

\noindent \textbf{2.1. Definition.} By a $d$-dimensional ($d\geq
2$) $A$-$D$-$E$ singularity we mean a singularity that is
analytically isomorphic to the germ at the origin of one of the
following hypersurfaces in $\mathbb{A}^{d+1}_{\mathbb{C}}$ (with
coordinates $(x_1,\ldots,x_{d+1})$):
\[\begin{array}{clccl}
(1)& x_1^{n+1}+x_2^2+x_3^2+\cdots+x_{d+1}^2=0 & & & (\text{type }
A_n,n\geq 1),\\ & & & & \\ (2)&
x_1^{n-1}+x_1x_2^2+x_3^2+\cdots+x_{d+1}^2=0 & & & (\text{type }
D_n,n\geq 4),\\ & & & & \\ (3)&
x_1^3+x_2^4+x_3^2+\cdots+x_{d+1}^2=0 & & & (\text{type } E_6),\\ & & & & \\
(4)& x_1^3+x_1x_2^3+x_3^2+\cdots+x_{d+1}^2=0 & & & (\text{type }
E_7),\\ & & & & \\ (5)& x_1^3+x_2^5+x_3^2+\cdots+x_{d+1}^2=0 & & &
(\text{type } E_8). \end{array}
\]

\noindent Some of their properties are listed in \cite[Remark
1.10]{DaisRoczen}.\\

\noindent \textbf{2.2.} We will now construct a log resolution for
these singularities by performing successive blow-ups, but we will
only do this for $d\geq 4$. The case $d=2$ is well known and the
construction in the three-dimensional case can be found in detail
in \cite[Section 2]{DaisRoczen}; in fact, our procedure is quite
analogous. The main differences are:
\begin{itemize}
\item[(1)] For $d\geq 4$, every blow-up adds just one component to
the exceptional locus, whereas you can get two planes intersecting
in a line as new exceptional divisors after a single blow-up in
the three-dimensional case (e.g. after the first blow-up in cases
$D$ and $E$).
\item[(2)] In the higher dimensional case, the analogue of this
line will be a singular line on the exceptional divisor, thus in
order to get a smooth normal crossings divisor one has to blow up
in such lines, which is not necessary for $d=3$.
\end{itemize}
An example will make this clear: blow up in the singular point of
the defining hypersurface in the $E_6$ case. For a suitable choice
of coordinates one finds $\{z_3^2+z_4^2=0\} \subset
\mathbb{P}^3_{\mathbb{C}}$ as equation of the exceptional locus
for $d=3$, and for $d\geq 4$, one finds
$\{z_3^2+z_4^2+\cdots+z_{d+1}^2=0\} \subset
\mathbb{P}^d_{\mathbb{C}}$ (this is irreducible, but the line
$\{z_3=\cdots=z_{d+1}=0\}$ is singular).

In what follows we use the same name for a divisor $D$ at the
moment of its creation as at all later stages (instead of speaking
of the strict transform of $D$). We work out the details for the
case of a $D_n$ singularity with even $n$ and we discuss the
results shortly in the other cases. We write $m$ for the number of
variables ($m\geq 5$) and use coordinates
$(x_1,\ldots,x_m)$ on $\mathbb{A}^m$.\\

\noindent \textbf{2.3.} \begin{tabular}{|c|} \hline \textbf{Case}
$\mathbf{A}$
\\ \hline
\end{tabular}\\

\noindent Consider the hypersurface
$X=\{x_1^{n+1}+x_2^2+\cdots+x_m^2=0\}\subset \mathbb{A}^m$ for
$m\geq 5$.\\

\noindent \underline{(1) $n$ odd, $n=2k-1$, with $k\geq
1$.}\\

Blowing up an $A_n$ singularity yields an $A_{n-2}$ singularity
(that lies on the exceptional locus) and nothing else happens.
Thus after $k$ point blow-ups we already have a log resolution.
The intersection diagram looks like

\begin{center}
\begin{picture}(70,10)
\put(2,2){\circle*{2}} \put(12,2){\circle*{2}}
\put(22,2){\circle*{2}} \put(52,2){\circle*{2}}
\put(62,2){\circle*{2}} \put(2,2){\line(1,0){26}}
\put(46,2){\line(1,0){16}} \put(0,6){$D_1$} \put(10,6){$D_2$}
\put(20,6){$D_3$} \put(49,6){$D_{k-1}$} \put(61,6){$D_k$}
\put(33,2){$\ldots$}
\end{picture}
\end{center}
where $D_i$ is created after the $i$-th blow-up. At the moment of
its creation, $D_i$ (for $i\in\{1,\ldots,k-1\}$) is isomorphic to
the singular quadric $\{x_2^2+\cdots+x_m^2=0\}$ in
$\mathbb{P}^{m-1}$, and its singular point is the center of the
next blow-up. The last divisor $D_k$ is isomorphic to the
nonsingular quadric in $\mathbb{P}^{m-1}$. In the end the
intersection of two exceptional divisors is isomorphic to a
nonsingular quadric in
$\mathbb{P}^{m-2}$.\\

\newpage

\noindent \underline{(2) $n$ even, $n=2k$, with $k\geq
1$.}\\

After $k$ point blow-ups the strict transform of $X$ is
nonsingular, but the last created divisor $D_k$ still has a
singular point, so we have to perform an extra blow-up (with
exceptional divisor $D_{k+1}$ isomorphic to $\mathbb{P}^{m-2}$).
As intersection diagram we find

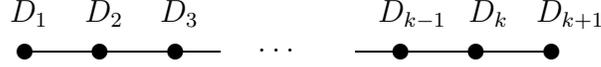
\begin{center}
\begin{picture}(81,10)
\put(3,2){\circle*{2}} \put(13,2){\circle*{2}}
\put(23,2){\circle*{2}} \put(53,2){\circle*{2}}
\put(63,2){\circle*{2}} \put(73,2){\circle*{2}}
\put(3,2){\line(1,0){26}} \put(47,2){\line(1,0){26}}
\put(1,6){$D_1$} \put(11,6){$D_2$} \put(21,6){$D_3$}
\put(50,6){$D_{k-1}$} \put(62,6){$D_k$} \put(71,6){$D_{k+1}$}
\put(34,2){$\ldots$}
\end{picture}
\end{center}
with all $D_i$ ($i\in \{1,\ldots,k\}$) isomorphic to the singular
quadric $\{x_2^2+\cdots+x_m^2=0\}$ in $\mathbb{P}^{m-1}$ at the
moment of their creation. Again, all intersections are isomorphic
to the
nonsingular quadric in $\mathbb{P}^{m-2}$.\\

\noindent \textbf{2.4.} \begin{tabular}{|c|} \hline \textbf{Case}
$\mathbf{D}$
\\ \hline
\end{tabular}\\

\noindent Now we study
$X=\{x_1^{n-1}+x_1x_2^2+x_3^2+\cdots+x_m^2=0\} \subset
\mathbb{A}^m$ for $m\geq 5$ and $n\geq 4$. Notice that you also
find singularities for $n=2$ and $n=3$, but they are analytically
isomorphic to two $A_1$ and one $A_3$ singularity respectively.\\

\noindent \underline{(1) $n$ even, $n=2k$, with $k\geq
2$.}\\

\emph{Step $1$}: We blow up $X$ in the origin. Take
$(x_1,\ldots,x_m)\times (z_1,\ldots,z_m)$ as coordinates on
$\mathbb{A}^m\times \mathbb{P}^{m-1}$. Consider the reducible
variety $X'$ in $\mathbb{A}^m\times \mathbb{P}^{m-1}$ given by the
equations
\[\left\{\begin{array}{lcc} x_1^{2k-1}+x_1x_2^2+x_3^2+\cdots+x_m^2=0
& & \\ x_iz_j=x_jz_i & & \forall\, i,j\in \{1,\ldots,m\}.
\end{array}\right.\]

In the open set $z_1\neq 0$, $X'$ is isomorphic to
$\{x_1^2(x_1^{2k-3}+x_1x_2^2+x_3^2+\cdots+x_m^2)=0\} \subset
\mathbb{A}^m$ by replacing $x_j$ by $x_1\frac{z_j}{z_1}$ and
renaming the affine coordinate $\frac{z_j}{z_1}$ as $x_j$ for
$j=2,\ldots,m$. The equation $x_1=0$ describes here the
exceptional locus, while the other equation gives us the strict
transform of $X$, in which we are interested. Their intersection
is the first exceptional divisor, we call it $D_1$. We can do the
same thing for any open set $z_i\neq 0$ and thus we can describe
$X'$ by the following set of equations:
\[\left\{\begin{array}{lcc} x_1^2(x_1^{2k-3}+x_1x_2^2+x_3^2+\cdots+x_m^2)=0 & & (1)\\
x_2^2(x_1^{2k-1}x_2^{2k-3}+x_1x_2+x_3^2+\cdots+x_m^2)=0 & & (2)\\
x_3^2(x_1^{2k-1}x_3^{2k-3}+x_1x_2^2x_3+1+x_4^2+\cdots+x_m^2)=0 & & (3)\\
\qquad\qquad \vdots & & \vdots\\
x_m^2(x_1^{2k-1}x_m^{2k-3}+x_1x_2^2x_m+x_3^2+\cdots+x_{m-1}^2+1)=0.
& & (m)
\end{array}\right.\]
One sees from this that globally $D_1 \cong
\{x_3^2+\cdots+x_m^2=0\}\subset \mathbb{P}^{m-1}$, which has a
singular line $\{x_3=\cdots=x_m=0\}$ (located in charts (1) and
(2)). Notice that for $k\geq 3$, we have a $D_{n-2}$ singularity
in chart (1) and a singularity that is analytically isomorphic to
an $A_1$ in the origin of chart (2). In the other charts both
$D_1$ and the strict transform of $X$ are nonsingular, so we have
no problems there. We will assume now that
$k\geq 4$ and we will see later what happens if $k=2,3$.\\

\emph{Step $2$}: Let us first get rid of the $A_1$ singularity.
Thus we blow up in the origin of chart (2). Since this blow-up is
an isomorphism outside this point, we preserve the other
coordinate charts and we replace chart (2) by the following
charts:
\[\left\{\begin{array}{lcc} x_1^4x_2^2(x_1^{4k-6}x_2^{2k-3}+x_2+x_3^2+\cdots+x_m^2)=0
 & & (2.1)\\
x_2^4(x_1^{2k-1}x_2^{4k-6}+x_1+x_3^2+\cdots+x_m^2)=0 & & (2.2)\\
x_2^2x_3^4(x_1^{2k-1}x_2^{2k-3}x_3^{4k-6}+x_1x_2+1+x_4^2+\cdots+x_m^2)=0 & & (2.3)\\
\qquad\qquad \vdots & & \vdots\\
x_2^2x_m^4(x_1^{2k-1}x_2^{2k-3}x_m^{4k-6}+x_1x_2+x_3^2+\cdots+x_{m-1}^2+1)=0.
& & (2.m)
\end{array}\right.\]
Now we see that the strict transform $\widetilde{X}$ of $X$ is
nonsingular in this part, but we still have the singular line on
$D_1$ (in charts (1) and (2.1) now). Our new exceptional divisor,
we call it $E_1$, is globally a nonsingular quadric in
$\mathbb{P}^{m-1}$.

We check immediately that $D_1$ and $E_1$ intersect transversally
outside the singular line of $D_1$: take a point
$P=(0,0,\alpha_3,\ldots,\alpha_m)$ on their intersection in chart
(2.1) for example (thus $\alpha_3^2+\cdots+\alpha_m^2=0$). We
assume that $P$ does not lie on the singular line on $D_1$ (so at
least one of the $\alpha_i$ is nonzero), since we will blow it up
later. The local ring $\mathcal{O}_{P,\widetilde{X}}$ is
isomorphic to
$\bigl(\frac{\mathbb{C}[x_1,\ldots,x_m]}{I}\bigr)_{m_P}$ with
$I=(x_1^{4k-6}x_2^{2k-3}+x_2+x_3^2+\cdots+x_m^2)$ and
$m_P=\frac{(x_1,x_2,x_3-\alpha_3,\ldots,x_m-\alpha_m)}{I}$. As a
$\mathbb{C}$-vector space, $\frac{m_P}{m_P^2}$ has dimension $m-1$
and is isomorphic to
$\frac{(x_1,x_2,x_3-\alpha_3,\ldots,x_m-\alpha_m)}{(x_1^2,x_1x_2,x_2^2,x_3^2-
2\alpha_3x_3+\alpha_3^2,\ldots)+I}$. It is generated by the set
$\{x_1,x_2,x_3-\alpha_3,\ldots,x_m-\alpha_m\}$ and the last $m-1$
generators are linearly dependent, since
\begin{center}
\begin{tabular}[t]{cccl}
\multicolumn{4}{l}{$x_2+2\alpha_3(x_3-\alpha_3)+\cdots+2\alpha_m(x_m-\alpha_m)$}\\
& &=&  $x_2+2\alpha_3x_3+\cdots+2\alpha_mx_m$\\
& &=& $x_1^{4k-6}x_2^{2k-3}+x_2+x_3^2+\cdots+x_m^2-(x_1^{4k-7}x_2^{2k-4})x_1x_2$\\
& & & $-(x_3^2-
2\alpha_3x_3+\alpha_3^2)-\cdots-(x_m^2- 2\alpha_mx_m+\alpha_m^2)$\\
& &=&0,
\end{tabular}
\end{center}
and thus $x_1$ and $x_2$ must be linearly independent. Hence $D_1$
and $E_1$ have normal crossings at
$(0,0,\alpha_3,\ldots,\alpha_m)$. Later on, we will not check the
normal crossings condition any more, it will be satisfied for all
divisors in the end.\\

\newpage

\emph{Step $3$}: We tackle the $D_{n-2}$ singularity in chart (1)
now. We blow up in its origin:
\[\left\{\begin{array}{lcc} x_1^4(x_1^{2k-5}+x_1x_2^2+x_3^2+\cdots+x_m^2)=0 & & (1.1)\\
x_1^2x_2^4(x_1^{2k-3}x_2^{2k-5}+x_1x_2+x_3^2+\cdots+x_m^2)=0 & & (1.2)\\
x_1^2x_3^4(x_1^{2k-3}x_3^{2k-5}+x_1x_2^2x_3+1+x_4^2+\cdots+x_m^2)=0 & & (1.3)\\
\qquad\qquad \vdots & & \vdots\\
x_1^2x_m^4(x_1^{2k-3}x_m^{2k-5}+x_1x_2^2x_m+x_3^2+\cdots+x_{m-1}^2+1)=0.
& & (1.m)
\end{array}\right.\]
It is no surprise that we find a $D_{n-4}$ singularity in the
origin of chart (1.1) and an $A_1$ in the origin of chart (1.2).
The newly created divisor, called $D_2$, intersects $D_1$ and has
a singular line in charts (1.1) and (1.2); the singular line of
$D_1$ from chart (1) is transferred to chart (1.2).\\

\emph{Step $4$}: We blow up in the origin of chart (1.2). The
singularity is resolved and the new divisor $E_2$ intersects both
$D_1$ and $D_2$:
\[\left\{\begin{array}{lcc} x_1^8x_2^4(x_1^{4k-10}x_2^{2k-5}+x_2+x_3^2+\cdots+x_m^2)=0
 & & (1.2.1)\\
x_1^2x_2^8(x_1^{2k-3}x_2^{4k-10}+x_1+x_3^2+\cdots+x_m^2)=0 & & (1.2.2)\\
x_1^2x_2^4x_3^8(x_1^{2k-3}x_2^{2k-5}x_3^{4k-10}+x_1x_2+1+x_4^2+\cdots+x_m^2)=0 & & (1.2.3)\\
\qquad\qquad \vdots & & \vdots\\
x_1^2x_2^4x_m^8(x_1^{2k-3}x_2^{2k-5}x_m^{4k-10}+x_1x_2+x_3^2+\cdots+x_{m-1}^2+1)=0.
& & (1.2.m)
\end{array}\right.\]
The singular lines on $D_1$ and $D_2$ are separated and go to
charts (1.2.2) and (1.2.1) respectively.\\

We continue in this way, performing alternate blow-ups in a $D_i$
and an $A_1$, until we
have to blow up in a $D_4$ singularity.\\

\emph{Step $n-3$}: We blow up in the origin of the chart
$x_1^{2k-4}(x_1^3+x_1x_2^2+x_3^2+\cdots+x_m^2)=0$.
\[\left\{\begin{array}{lcc} x_1^{2k-2}(x_1+x_1x_2^2+x_3^2+\cdots+x_m^2)=0 & & (1')\\
x_1^{2k-4}x_2^{2k-2}(x_1^3x_2+x_1x_2+x_3^2+\cdots+x_m^2)=0 & & (2')\\
x_1^{2k-4}x_3^{2k-2}(x_1^3x_3+x_1x_2^2x_3+1+x_4^2+\cdots+x_m^2)=0 & & (3')\\
\qquad\qquad \vdots & & \vdots\\
x_1^{2k-4}x_m^{2k-2}(x_1^3x_m+x_1x_2^2x_m+x_3^2+\cdots+x_{m-1}^2+1)=0.
& & (m')
\end{array}\right.\]
In fact $(j')$ stands here
for\begin{tabular}[t]{c}($\underbrace{1.1\ldots 1}.j$)\\ $k-2$
times
\end{tabular}. We get three singular points, all an\-a\-ly\-ti\-cal\-ly
isomorphic to an $A_1$ singularity. Both present divisors (we call
them of course $D_{k-2}$ and $D_{k-1}$) have a singular line and
in fact all the singular points lie on the singular line of
$D_{k-1}$. One of the singular points, the origin of chart $(2')$,
lies on the intersection of $D_{k-2}$ and $D_{k-1}$. Note that the
singular points $(0,i,0,\ldots,0)$ and $(0,-i,0,\ldots,0)$ of
chart $(1')$ correspond to the points $(-i,0,\ldots,0)$ and
$(i,0,\ldots,0)$ of chart $(2')$
respectively. \\

\emph{Step $n-2$}: We deal with the origin of chart $(2')$ first.
Blowing it up yields a divisor $E_{k-1}$ that intersects $D_{k-1}$
and $D_{k-2}$:
\[\left\{\begin{array}{lcc} x_1^{4k-4}x_2^{2k-2}(x_1^2x_2+x_2+x_3^2+\cdots+x_m^2)=0 &
& (2'.1)\\ x_1^{2k-4}x_2^{4k-4}(x_1^3x_2^2+x_1+x_3^2+\cdots+x_m^2)=0 & & (2'.2)\\
x_1^{2k-4}x_2^{2k-2}x_3^{4k-4}(x_1^3x_2x_3^2+x_1x_2+1+x_4^2+\cdots+x_m^2)=0 & & (2'.3)\\
\qquad\qquad \vdots & & \vdots\\
x_1^{2k-4}x_2^{2k-2}x_m^{4k-4}(x_1^3x_2x_m^2+x_1x_2+x_3^2+\cdots+x_{m-1}^2+1)=0.
& & (2'.m)
\end{array}\right.\]
The other two singularities lie in charts $(1')$ and $(2'.1)$. The singular lines on
$D_{k-2}$ and $D_{k-1}$ get separated and go to charts $(2'.2)$ and $(2'.1)$, respectively.\\

\emph{Step $n-1$}: After a coordinate transformation the equation
of chart $(1')$ becomes
$x_1^{2k-2}(x_1x_2(x_2+2i)+x_3^2+\cdots+x_m^2=0$. To put the same
point in the origin, we have to change the equation of chart
$(2'.1)$ to
$(x_1-i)^{4k-4}x_2^{2k-2}(x_1x_2(x_1-2i)+x_3^2+\cdots+x_m^2)=0$
for example. In this step we blow up both charts in the origin and
we call the new divisor $F_1$:
\begin{equation*}
\begin{split}
&\left\{\begin{array}{lcc} x_1^{2k}(x_2(x_1x_2+2i)+x_3^2+\cdots+x_m^2)=0 & & (1'.1)\\
x_1^{2k-2}x_2^{2k}(x_1(x_2+2i)+x_3^2+\cdots+x_m^2)=0 & & (1'.2)\\
x_1^{2k-2}x_3^{2k}(x_1x_2(x_2x_3+2i)+1+x_4^2+\cdots+x_m^2)=0 & & (1'.3)\\
\qquad\qquad \vdots & & \vdots\\
x_1^{2k-2}x_m^{2k}(x_1x_2(x_2x_m+2i)+x_3^2+\cdots+x_{m-1}^2+1)=0 &
& (1'.m)
\end{array}\right. \quad \text{ and } \\
&\quad \left\{\begin{array}{lcc}
x_1^{2k}(x_1-i)^{4k-4}x_2^{2k-2}(x_2(x_1-2i)+x_3^2+\cdots+x_m^2)=0
& & (2'.1.1)\\ (x_1x_2-i)^{4k-4}x_2^{2k}(x_1(x_1x_2-2i)+x_3^2+\cdots+x_m^2)=0 & & (2'.1.2)\\
(x_1x_3-i)^{4k-4}x_2^{2k-2}x_3^{2k}(x_1x_2(x_1x_3-2i)+1+\cdots+x_m^2)=0 & & (2'.1.3)\\
\qquad\qquad \vdots & & \vdots\\
(x_1x_m-i)^{4k-4}x_2^{2k-2}x_m^{2k}(x_1x_2(x_1x_m-2i)+x_3^2+\cdots+1)=0.
& & (2'.1.m)
\end{array}\right. \end{split}
\end{equation*}
The last singular point and the singular line on $D_{k-1}$ are now
in charts $(1'.2)$ and $(2'.1.1)$.\\

\emph{Step $n$}: Before blowing up the final singular point, we
first do a coordinate transformation in chart $(1'.2)$ to get the
equation $x_1^{2k-2}(x_2-2i)^{2k}(x_1x_2+x_3^2+\cdots+x_m^2)=0$
and in chart $(2'.1.1)$ to get
$(x_1+2i)^{2k}(x_1+i)^{4k-4}x_2^{2k-2}(x_1x_2+x_3^2+\cdots+x_m^2)=0$.
The new exceptional divisor is called $F_2$.
\begin{equation*}
\begin{split}
&\left\{\begin{array}{lcc} x_1^{2k}(x_1x_2-2i)^{2k}(x_2+x_3^2+\cdots+x_m^2)=0 & & (1'.2.1)\\
x_1^{2k-2}(x_2-2i)^{2k}x_2^{2k}(x_1+x_3^2+\cdots+x_m^2)=0 & & (1'.2.2)\\
x_1^{2k-2}(x_2x_3-2i)^{2k}x_3^{2k}(x_1x_2+1+x_4^2+\cdots+x_m^2)=0 & & (1'.2.3)\\
\qquad\qquad \vdots & & \vdots\\
x_1^{2k-2}(x_2x_m-2i)^{2k}x_m^{2k}(x_1x_2+x_3^2+\cdots+x_{m-1}^2+1)=0
& & (1'.2.m)
\end{array}\right. \quad \text{ and } \\
&\ \ \left\{\begin{array}{lcc}
x_1^{2k}(x_1+2i)^{2k}(x_1+i)^{4k-4}x_2^{2k-2}(x_2+x_3^2+\cdots+x_m^2)=0
& & (2'.1.1.1)\\
(x_1x_2+2i)^{2k}(x_1x_2+i)^{4k-4}x_2^{2k}(x_1+x_3^2+\cdots+x_m^2)=0
& &
(2'.1.1.2)\\
(x_1x_3+2i)^{2k}(x_1x_3+i)^{4k-4}x_2^{2k-2}x_3^{2k}(x_1x_2+1+\cdots+x_m^2)=0
& &
(2'.1.1.3)\\ \qquad\qquad \vdots & & \vdots\\
(x_1x_m+2i)^{2k}(x_1x_m+i)^{4k-4}x_2^{2k-2}x_m^{2k}(x_1x_2+x_3^2+\cdots+1)=0.
& & (2'.1.1.m)
\end{array}\right. \end{split}
\end{equation*}
The singular line on $D_{k-1}$ is moved to charts $(1'.2.2)$ and
$(2'.1.1.1)$.\\

In the next $k-1$ steps we blow up in the singular lines on the
divisors $D_i$. This gives rise to new exceptional divisors which
will be denoted by $G_i$. After $k-1$ steps we finally have a log
resolution; we will perform steps $n+1$ and $n+k-1$ explicitly.\\

\emph{Step $n+1$}: To cover the singular line on $D_1$ completely,
we have to perform the blow-up in charts (2.1) and (1.2.2). In
chart (2.1) we have to blow up the variety
$Y=\{x_1^4x_2^2(x_1^{4k-6}x_2^{2k-3}+x_2+x_3^2+\cdots+x_m^2)=0\}\subset
\mathbb{A}^m$ in the line $\{x_2=\cdots=x_m=0\}$. The strict
transform of $Y$ and the exceptional locus form a reducible
variety in $\mathbb{A}^m\times \mathbb{P}^{m-2}$, given by the
equations
\[\left\{\begin{array}{lcc} x_1^4x_2^2(x_1^{4k-6}x_2^{2k-3}+x_2+x_3^2+\cdots+x_m^2)=0
& & \\ x_iz_j=x_jz_i & & \forall\, i,j\in \{2,\ldots,m\},
\end{array}\right.\]
where $(z_2,\ldots,z_m)$ are homogenous coordinates on
$\mathbb{P}^{m-2}$. As for a point blow-up, we can replace $x_j$
by $x_i\frac{z_j}{z_i}$ in the open set $z_i\neq 0$ and rename
$\frac{z_j}{z_i}$ as $x_j$. Hence we get the following equations
for $Y'$:
\[\left\{\begin{array}{lcc} x_1^4x_2^3(x_1^{4k-6}x_2^{2k-4}+1+x_2x_3^2+\cdots+x_2x_m^2)=0
 & & (2.1.2)\\
x_1^4x_2^2x_3^3(x_1^{4k-6}x_2^{2k-3}x_3^{2k-4}+x_2+x_3+x_3x_4^2+\cdots+x_3x_m^2)=0 & &
(2.1.3)\\ \qquad\qquad \vdots & & \vdots\\
x_1^4x_2^2x_m^3(x_1^{4k-6}x_2^{2k-3}x_m^{2k-4}+x_2+x_3^2x_m+\cdots+x_{m-1}^2x_m+x_m)=0.
& & (2.1.m)
\end{array}\right.\]

The equations after blowing up in $\{x_1=x_3=\cdots=x_m=0\}$ in
chart (1.2.2) are:
\[\left\{\begin{array}{lcc} x_1^3x_2^8(x_1^{2k-4}x_2^{4k-10}+1+x_1x_3^2+\cdots+x_1x_m^2)=0
 & & (1.2.2.1)\\
x_1^2x_2^8x_3^3(x_1^{2k-3}x_2^{4k-10}x_3^{2k-4}+x_1+x_3+x_3x_4^2+\cdots+x_3x_m^2)=0
& &
(1.2.2.3)\\ \qquad\qquad \vdots & & \vdots\\
x_1^2x_2^8x_m^3(x_1^{2k-3}x_2^{4k-10}x_m^{2k-4}+x_1+x_3^2x_m+\cdots+x_{m-1}^2x_m+x_m)=0.
& & (1.2.2.m)
\end{array}\right.\]

\emph{Step $n+k-1$}: Here we have to consider charts $(1'.2.2)$
and $(2'.1.1.1)$ in which $D_{k-1}$ still has a singular line with
equations $\{x_1=x_3=\cdots =x_m=0\}$ and $\{x_2=x_3=\cdots
=x_m=0\}$, respectively. Blowing it up yields
\begin{equation*}
\begin{split}
&\left\{\begin{array}{lcc}
x_1^{2k-1}(x_2-2i)^{2k}x_2^{2k}(1+x_1x_3^2+\cdots+x_1x_m^2)=0 & &
(1'.2.2.1)\\
x_1^{2k-2}(x_2-2i)^{2k}x_2^{2k}x_3^{2k-1}(x_1+x_3+\cdots+x_3x_m^2)=0
& &
(1'.2.2.3)\\ \qquad\qquad \vdots & & \vdots\\
x_1^{2k-2}(x_2-2i)^{2k}x_2^{2k}x_m^{2k-1}(x_1+x_3^2x_m+\cdots+x_m)=0
& & (1'.2.2.m)
\end{array}\right. \ \ \text{ and } \\
&\ \ \left\{\begin{array}{lc}
x_1^{2k}(x_1+2i)^{2k}(x_1+i)^{4k-4}x_2^{2k-1}(1+x_2x_3^2+\cdots+x_2x_m^2)=0
&  (2'.1.1.1.2)\\
x_1^{2k}(x_1+2i)^{2k}(x_1+i)^{4k-4}x_2^{2k-2}x_3^{2k-1}(x_2+x_3+\cdots+x_3x_m^2)=0
&
(2'.1.1.1.3)\\ \qquad\qquad \vdots  & \vdots\\
x_1^{2k}(x_1+2i)^{2k}(x_1+i)^{4k-4}x_2^{2k-2}x_m^{2k-1}(x_2+x_3^2x_m+\cdots
+x_m)=0.  & (2'.1.1.1.m)
\end{array}\right. \end{split}
\end{equation*}\vspace{2mm}

From these calculations, we can deduce the intersection diagram.
We leave it to the reader to check the details. It can be easily
seen that the same diagram is valid for $k=2,3$.

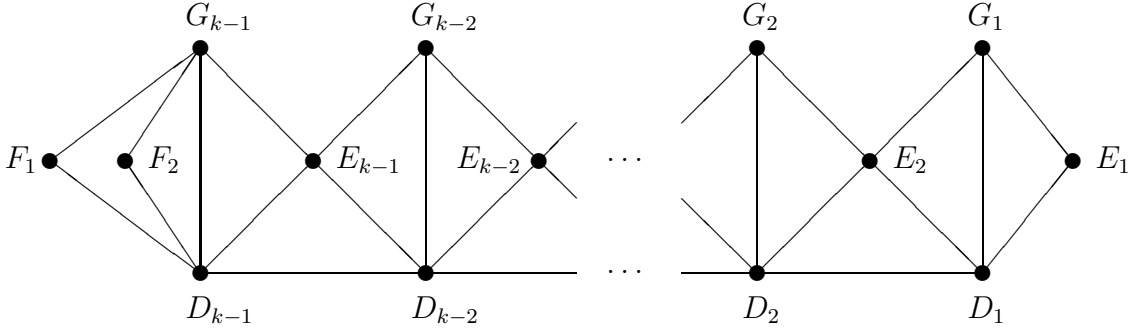
\begin{center}
\begin{picture}(160,45)
\put(27,7){\circle*{2}} \put(57,7){\circle*{2}}
\put(101,7){\circle*{2}} \put(131,7){\circle*{2}}
\put(27,37){\circle*{2}} \put(57,37){\circle*{2}}
\put(101,37){\circle*{2}} \put(131,37){\circle*{2}}
\put(42,22){\circle*{2}} \put(72,22){\circle*{2}}
\put(116,22){\circle*{2}} \put(143,22){\circle*{2}}
\put(7,22){\circle*{2}} \put(17,22){\circle*{2}}
\put(27,7){\line(1,0){50}} \put(91,7){\line(1,0){40}}
\put(27,7){\line(0,1){30}} \put(57,7){\line(0,1){30}}
\put(101,7){\line(0,1){30}} \put(131,7){\line(0,1){30}}
\put(27,7){\line(1,1){30}} \put(57,7){\line(1,1){20}}
\put(101,7){\line(1,1){30}} \put(131,7){\line(4,5){12}}
\put(101,37){\line(1,-1){30}} \put(131,37){\line(4,-5){12}}
\put(27,37){\line(1,-1){30}} \put(57,37){\line(1,-1){20}}
\put(101,7){\line(-1,1){10}} \put(101,37){\line(-1,-1){10}}
\put(27,7){\line(-2,3){10}} \put(27,7){\line(-4,3){20}}
\put(17,22){\line(2,3){10}} \put(7,22){\line(4,3){20}}
\put(25,40){$G_{k-1}$} \put(55,40){$G_{k-2}$} \put(99,40){$G_2$}
\put(129,40){$G_{1}$} \put(25,1){$D_{k-1}$} \put(55,1){$D_{k-2}$}
\put(99,1){$D_{2}$} \put(129,1){$D_{1}$} \put(1,21){$F_1$}
\put(20,21){$F_2$} \put(45,21){$E_{k-1}$} \put(61,21){$E_{k-2}$}
\put(119,21){$E_2$} \put(146,21){$E_1$} \put(81,22){$\ldots$}
\put(81,7){$\ldots$}
\end{picture}
\end{center}

$\phantom{some place}$

\noindent \underline{(2) $n$ odd, $n=2k+1$, with $k\geq
2$.}\\

The first $2k-4$ steps are completely analogous to the case where
$n$ is even. Now we end up with the equation
$x_1^{2k-4}(x_1^4+x_1x_2^2+x_3^2+\cdots+x_m^2)$ which has a $D_5$
singularity in the origin. Blowing this up gives one $A_3$
singularity on the new divisor $D_{k-1}$ (the equation of the
first chart is $x_1^{2k-2}(x_1^2+x_1x_2^2+x_3^2+\cdots+x_m^2=0)$).
We already know that this can be resolved by two consecutive
blow-ups, creating divisors $F_1$ and $F_2$. Afterwards, the
singular lines on the $D_i$ must be blown up. Explicit
calculations will lead to the following intersection diagram:

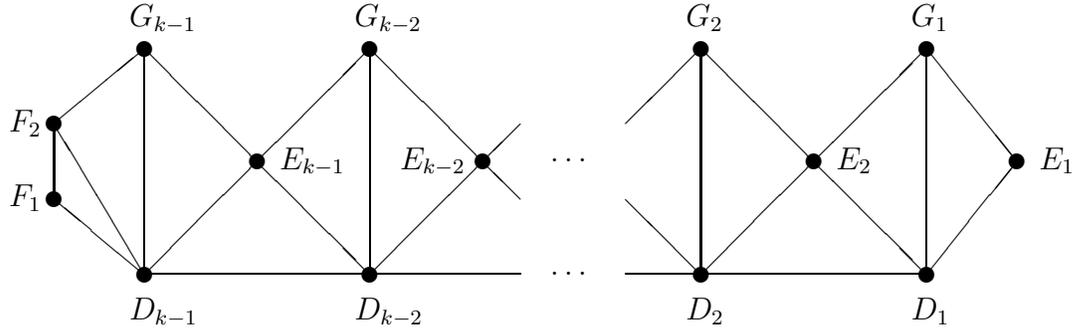
\begin{center}
\begin{picture}(160,45)
\put(19,7){\circle*{2}} \put(49,7){\circle*{2}}
\put(93,7){\circle*{2}} \put(123,7){\circle*{2}}
\put(19,37){\circle*{2}} \put(49,37){\circle*{2}}
\put(93,37){\circle*{2}} \put(123,37){\circle*{2}}
\put(34,22){\circle*{2}} \put(64,22){\circle*{2}}
\put(108,22){\circle*{2}} \put(135,22){\circle*{2}}
\put(7,17){\circle*{2}} \put(7,27){\circle*{2}}
\put(19,7){\line(1,0){50}} \put(83,7){\line(1,0){40}}
\put(19,7){\line(0,1){30}} \put(49,7){\line(0,1){30}}
\put(93,7){\line(0,1){30}} \put(123,7){\line(0,1){30}}
\put(19,7){\line(1,1){30}} \put(49,7){\line(1,1){20}}
\put(93,7){\line(1,1){30}} \put(123,7){\line(4,5){12}}
\put(93,37){\line(1,-1){30}} \put(123,37){\line(4,-5){12}}
\put(19,37){\line(1,-1){30}} \put(49,37){\line(1,-1){20}}
\put(93,7){\line(-1,1){10}} \put(93,37){\line(-1,-1){10}}
\put(19,7){\line(-6,5){12}} \put(19,7){\line(-3,5){12}}
\put(7,27){\line(6,5){12}} \put(7,17){\line(0,1){10}}
\put(17,40){$G_{k-1}$} \put(47,40){$G_{k-2}$} \put(91,40){$G_2$}
\put(121,40){$G_{1}$} \put(17,1){$D_{k-1}$} \put(47,1){$D_{k-2}$}
\put(91,1){$D_{2}$} \put(121,1){$D_{1}$} \put(1,16){$F_1$}
\put(1,26){$F_2$} \put(37,21){$E_{k-1}$} \put(53,21){$E_{k-2}$}
\put(111,21){$E_2$} \put(138,21){$E_1$} \put(73,22){$\ldots$}
\put(73,7){$\ldots$}
\end{picture}
\end{center}

$\phantom{some place}$

\noindent \textbf{2.5.} \begin{tabular}{|c|} \hline \textbf{Case}
$\mathbf{E_6}$
\\ \hline
\end{tabular}\\

\noindent After blowing up in the origin we get an $A_5$
singularity and a singular line on the first exceptional divisor
$D_1$. To resolve the $A_5$ singularity we need three more point
blow-ups (creating $D_2,D_3$ and $D_4$) and in the end we blow up
in the singular line (giving rise to a divisor $D_5$). We find
as intersection graph:\\
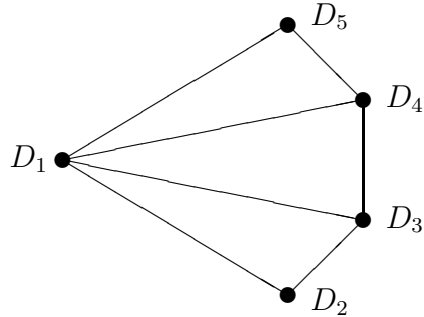
\begin{center}
\begin{picture}(60,38)
\put(10,20){\circle*{2}} \put(40,2){\circle*{2}}
\put(40,38){\circle*{2}} \put(50,12){\circle*{2}}
\put(50,28){\circle*{2}} \put(10,20){\line(5,-3){30}}
\put(10,20){\line(5,3){30}} \put(10,20){\line(5,1){40}}
\put(10,20){\line(5,-1){40}} \put(40,2){\line(1,1){10}}
\put(40,38){\line(1,-1){10}} \put(50,12){\line(0,1){16}}
\put(3,19){$D_1$} \put(43,0){$D_2$} \put(53,11){$D_3$}
\put(53,27){$D_4$} \put(43,38){$D_5$}
\end{picture}
\end{center}

\noindent \textbf{2.6.} \begin{tabular}{|c|} \hline \textbf{Cases}
$\mathbf{E_7}$ \textbf{and} $\mathbf{E_8}$
\\ \hline
\end{tabular}\\

\noindent An $E_7$ becomes a $D_6$ after one step and calculating
the intersections gives the following diagram
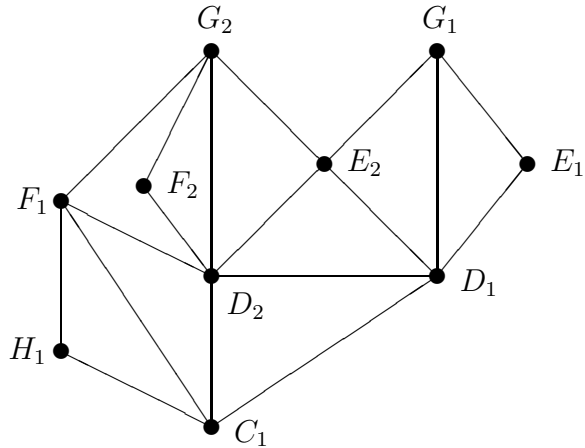
\begin{center}
\begin{picture}(80,57)
\put(10,32){\circle*{2}} \put(21,34){\circle*{2}}
\put(30,22){\circle*{2}} \put(30,52){\circle*{2}}
\put(45,37){\circle*{2}} \put(60,22){\circle*{2}}
\put(60,52){\circle*{2}} \put(72,37){\circle*{2}}
\put(30,2){\circle*{2}} \put(10,12){\circle*{2}}
\put(30,2){\line(-2,1){20}} \put(30,2){\line(-2,3){20}}
\put(30,2){\line(0,1){50}} \put(30,2){\line(3,2){30}}
\put(10,12){\line(0,1){20}} \put(10,32){\line(2,-1){20}}
\put(10,32){\line(1,1){20}} \put(21,34){\line(1,2){9}}
\put(21,34){\line(3,-4){9}} \put(30,22){\line(1,1){30}}
\put(30,52){\line(1,-1){30}} \put(30,22){\line(1,0){30}}
\put(60,22){\line(0,1){30}} \put(60,22){\line(4,5){12}}
\put(60,52){\line(4,-5){12}} \put(33,0){$C_1$} \put(63,20){$D_1$}
\put(32,17){$D_2$} \put(75,36){$E_1$} \put(48,36){$E_2$}
\put(4,31){$F_1$} \put(24,33){$F_2$} \put(28,55){$G_2$}
\put(58,55){$G_1$} \put(3,11){$H_1$}
\end{picture}
\end{center}

\noindent where $C_1$ is the very first exceptional divisor and
where $H_1$ arises after blowing up the singular line on $C_1$.
The other divisors come from the $D_6$ singularity. Notice the
difference between $F_1$ and $F_2$. It is easy to see that an
$E_8$ singularity passes to an $E_7$ after one blow-up, with again
a singular line on the first exceptional divisor $B_1$. We denote
the divisor that appears after blowing up in this singular line by
$I_1$ and we find the following intersection graph:
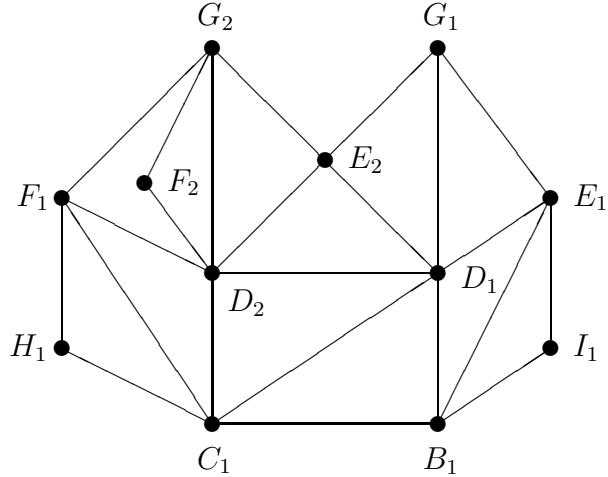
\begin{center}
\begin{picture}(83,64)
\put(10,37){\circle*{2}} \put(21,39){\circle*{2}}
\put(30,27){\circle*{2}} \put(30,57){\circle*{2}}
\put(45,42){\circle*{2}} \put(60,27){\circle*{2}}
\put(60,57){\circle*{2}} \put(75,37){\circle*{2}}
\put(30,7){\circle*{2}} \put(10,17){\circle*{2}}
\put(60,7){\circle*{2}} \put(75,17){\circle*{2}}
\put(30,7){\line(-2,1){20}} \put(30,7){\line(-2,3){20}}
\put(30,7){\line(0,1){50}} \put(30,7){\line(3,2){30}}
\put(10,17){\line(0,1){20}} \put(10,37){\line(2,-1){20}}
\put(10,37){\line(1,1){20}} \put(21,39){\line(1,2){9}}
\put(21,39){\line(3,-4){9}} \put(30,27){\line(1,1){30}}
\put(30,57){\line(1,-1){30}} \put(30,27){\line(1,0){30}}
\put(60,27){\line(0,1){30}} \put(60,27){\line(3,2){15}}
\put(60,57){\line(3,-4){15}} \put(30,7){\line(1,0){30}}
\put(60,7){\line(0,1){20}} \put(60,7){\line(1,2){15}}
\put(60,7){\line(3,2){15}} \put(75,17){\line(0,1){20}}
\put(28,1){$C_1$} \put(63,25){$D_1$} \put(32,22){$D_2$}
\put(78,36){$E_1$} \put(48,41){$E_2$} \put(4,36){$F_1$}
\put(24,38){$F_2$} \put(28,60){$G_2$} \put(58,60){$G_1$}
\put(3,16){$H_1$} \put(78,16){$I_1$} \put(58,1){$B_1$}
\end{picture}
\end{center}

$\phantom{some place}$

\section{The Hodge-Deligne polynomials of the pieces of the exceptional locus}

\noindent \textbf{3.1.} Denote by $a_r,b_r,c_r$ ($r\geq 2$) the
Hodge-Deligne polynomials of
\begin{itemize}
\item
$\{x_1^2+\cdots+x_r^2=0\}\subset \mathbb{P}^{r+1}_{\mathbb{C}}$,
\item $\{x_1^2+\cdots+x_r^2=0\}\subset \mathbb{P}^r_{\mathbb{C}}$,
\item $\{x_1^2+\cdots+x_r^2=0\}\subset
\mathbb{P}^{r-1}_{\mathbb{C}}$,
\end{itemize}
respectively, where $\mathbb{P}^s$ gets coordinates
$(x_1,\ldots,x_{s+1})$. We will be able to express all the needed
Hodge-Deligne polynomials in terms of $a_r,b_r$ and $c_r$, and
these last expressions are well known. For completeness we include
their computation in the following lemma. \emph{From now on, we
will write $w$ as abbreviation of $uv$}.

\newpage

\noindent \textbf{Lemma.} \textsl{The formulae for $a_r,b_r$ and
$c_r$ are given in the following table}:
\begin{center}
\begin{tabular}[t]{c|c|c}
 & $r$ \textsl{even} & $r$  \textsl{odd} \\
\hline & & \\ $a_r$ & $\frac{w^{r+1}-1}{w-1}+w^{\frac{r}{2}+1} $ &
$\frac{w^{r+1}-1}{w-1} $\\ & & \\
\hline & &  \\ $b_r$ & $\frac{w^r-1}{w-1}+w^{\frac{r}{2}} $ & $\frac{w^r-1}{w-1} $\\ & & \\
\hline & & \\ $c_r$ & $\frac{w^{r-1}-1}{w-1}+w^{\frac{r}{2}-1} $ &
$\frac{w^{r-1}-1}{w-1} $\\ & &
\end{tabular}
\end{center}

\noindent \textbf{Proof:} Denote by $d_r$ the Hodge-Deligne polynomial of
$\{x_1^2+\cdots+x_r^2+1=0\}\subset \mathbb{A}^r$. First we compute
$d_r$ by induction on $r$. Since $d_2$ is the Hodge-Deligne
polynomial of a conic with two points at infinity, it equals
$w-1$. The variety $\{x_1^2+x_2^2+x_3^2+1=0\} \subset
\mathbb{A}^3$ can be regarded as $\mathbb{P}^1\times \mathbb{P}^1$
minus a conic and thus $d_3=(w+1)^2-(w+1)=w^2+w$. For $r\geq 4$ we
use the isomorphism $\{x_1^2+\cdots +x_r^2+1=0\}\cong
\{x_1x_2+x_3^2+\cdots+x_r^2+1=0\}$. If $x_1=0$ in this last
equation, then the contribution to $d_r$ is $wd_{r-2}$ and if
$x_1\neq 0$, then it is $(w-1)w^{r-2}$, so we have the recursion
formula $d_r=wd_{r-2}+(w-1)w^{r-2}$. From this it follows that
$d_r=w^{r-1}-w^{\frac{r}{2}-1}$ if $r$ is even and
$d_r=w^{r-1}+w^{\frac{r-1}{2}}$ if $r$ is odd.

For $a_2$ we find $2w^2+w+1$ and we have the recursion formula
$a_r=a_{r-1}+w^2d_{r-1}$ for $r\geq 3$. The formulae for $b_r$ and
$c_r$ can be deduced similarly. \hfill $\blacksquare$\\

\noindent \textbf{3.2.} For the remainder of this section, we will
calculate the Hodge-Deligne polynomials of the pieces
$D_{J}^{\circ}$ (see the definition of the stringy $E$-function).
Since we are mainly interested in the contribution of the singular
point (by which we mean $E_{st}(X) -
H(D_{\emptyset}^{\circ})=E_{st}(X)-H(X\setminus \{0\})$, where $X$ is a
defining variety of an $A$-$D$-$E$ singularity), we will do this
for $J\neq \emptyset$.

We remark here the following. In the defining formula of the
stringy $E$-function we need the Hodge-Deligne polynomials of the
$D_J^{\circ}$ at the end of the resolution process. Notice however
that we can compute them immediately after they are created, since
a blow-up is an isomorphism outside its center. So we just have to
subtract contributions of intersections with previously created
divisors and already present centers of future blow-ups from the
global Hodge-Deligne polynomial in the right way.

The case of an $A$-$D$-$E$ surface singularity is well known and
for threefold singularities we refer again to \cite{DaisRoczen},
so we consider here the higher dimensional case. Parallel to the
previous section, we will work out the details for the case $D_n$,
$n$ even, and state the results in the other cases. We use the
same notations as in the previous
section.\\

\newpage

\noindent \textbf{3.3.} \begin{tabular}{|c|} \hline \textbf{Case}
$\mathbf{A}$
\\ \hline
\end{tabular}\\

\noindent From the description in (2.3), one gets the
following:\\

\noindent \underline{(1) $n$ odd}

\[\begin{array}{lccr}
H(D_1^{\circ})=b_{m-1}-1 & & & \\
H(D_i^{\circ})=b_{m-1}-c_{m-1}-1 & & &(i=2,\ldots,k-1)\\
H(D_k^{\circ})=c_m-c_{m-1} & & &\\
H(D_i\cap D_{i+1})=c_{m-1} & & &(i=1,\ldots,k-1)\\
\end{array}  \]

\noindent \underline{(2) $n$ even}

\[\begin{array}{lccr}
H(D_1^{\circ})=b_{m-1}-1 & & & \\
H(D_i^{\circ})=b_{m-1}-c_{m-1}-1 & & &(i=2,\ldots,k)\\
H(D_{k+1}^{\circ})=w^{m-2}+\cdots+1-c_{m-1} & & &\\
H(D_i\cap D_{i+1})=c_{m-1} & & &(i=1,\ldots,k)\\
\end{array}  \]

$\phantom{some place}$

\noindent \textbf{3.4.} \begin{tabular}{|c|} \hline \textbf{Case}
$\mathbf{D}$
\\ \hline
\end{tabular}\\

\noindent \underline{(1) $n$ even}\\

All the needed information can be read off from the equations in
(2.4). We follow the same steps.\\

\emph{Step $1$}: The first exceptional divisor is globally
isomorphic to $\{x_3^2+\cdots+x_m^2=0\}\subset \mathbb{P}^{m-1}$,
which has a singular line that contains the two singular points of
the surrounding variety. Hence $H(D_1^{\circ})=a_{m-2}-(w+1)$.\\

\emph{Step $2$}: One sees that $E_1$ is a nonsingular quadric in
$\mathbb{P}^{m-1}$ that intersects $D_1$ in
$\{x_3^2+\cdots+x_m^2=0\}\subset\mathbb{P}^{m-2}$, for a suitable
choice of coordinates. Thus $H(E_1^{\circ})=c_{m}-b_{m-2}$. The
intersection of $D_1$ and $E_1$ contains one point of the singular
line on $D_1$ and hence $H((D_1\cap E_1)^{\circ})=b_{m-2}-1$.\\

\emph{Step $3$}: Analogous to step 1 one finds that $D_2$ is
isomorphic to $\{x_3^2+\cdots+x_m^2=0\}\subset \mathbb{P}^{m-1}$,
with a singular line that contains two singular points of the
surrounding variety. Now $D_2$ intersects $D_1$ in
$\{x_3^2+\cdots+x_m^2=0\}\subset \mathbb{P}^{m-2}$. This
intersection has exactly one point (the origin of coordinate chart
(1.2)) in common with the singular lines on $D_2$ and $D_1$. The
conclusion is that $H(D_2^{\circ})=a_{m-2}-(w+1)-b_{m-2}+1$ and
$H((D_1\cap D_2)^{\circ})=b_{m-2}-1$.\\

\emph{Step $4$}: For $H(E_2^{\circ})$ we find
$c_m-2b_{m-2}+c_{m-2}$, where $2b_{m-2}$ comes from the
intersections with $D_1$ and $D_2$ and $c_{m-2}$ from the
intersection with $D_1\cap D_2$. We also have that $H((D_1\cap
E_2)^{\circ})=H((D_2\cap E_2)^{\circ})=b_{m-2}-c_{m-2}-1$, where
the $-1$ comes from a point on the singular lines on the $D_i$.
Finally $H(D_1\cap D_2\cap E_2)=c_{m-2}$.\\

Analogously, for all $i$ from 3 to $k-2$, we have
$H(D_i^{\circ})=a_{m-2}-(w+1)-b_{m-2}+1$, $H((D_{i-1}\cap
D_{i})^{\circ})=b_{m-2}-1$, $H(E_i^{\circ})=c_m-2b_{m-2}+c_{m-2}$,
$H((D_{i-1}\cap E_{i})^{\circ})=H((D_i\cap
E_{i})^{\circ})=b_{m-2}-c_{m-2}-1$ and
$H(D_{i-1}\cap D_{i}\cap E_i)=c_{m-2}$.\\

\emph{Step $n-3$}: In this step three singular points are created,
but since they are all on the singular line on $D_{k-1}$, we still
find $H(D_{k-1}^{\circ})=a_{m-2}-(w+1)-b_{m-2}+1$ and
$H((D_{k-2}\cap D_{k-1})^{\circ})=b_{m-2}-1$.\\

\emph{Step $n-2$}: Again nothing special happens:
$H(E_{k-1}^{\circ})=c_m-2b_{m-2}+c_{m-2}$, $H((D_{k-2}\cap
E_{k-1})^{\circ})=H((D_{k-1}\cap
E_{k-1})^{\circ})=b_{m-2}-c_{m-2}-1$ and $H(D_{k-2}\cap
D_{k-1}\cap E_{k-1})=c_{m-2}$.\\

\emph{Step $n-1$ and step $n$}: Both $F_1$ and $F_2$ are
nonsingular quadrics in $\mathbb{P}^{m-1}$ and their intersection
with $D_{k-1}$ is $\{x_3^2+\cdots+x_m^2=0\}\subset
\mathbb{P}^{m-2}$, which has one point in common with the singular
line on $D_{k-1}$. Thus
$H(F_1^{\circ})=H(F_2^{\circ})=c_m-b_{m-2}$ and $H((D_{k-1}\cap
F_1)^{\circ})=H((D_{k-1}\cap F_2)^{\circ})=b_{m-2}-1$.\\

\emph{Step $n+1$}: The singular line on $D_1$ is except for the
origin of coordinate chart (2.1) covered by chart (1.2.2). But
after the blow-up, exactly the intersection of $E_1$ and $G_1$
lies above the origin of chart (2.1). Thus to calculate
$H(G_1^{\circ})$, it suffices to consider only charts (1.2.2.1) to
$(1.2.2.m)$. In chart (1.2.2.3) $G_1$ is just isomorphic to
$\mathbb{A}^{m-2}$. The piece of $G_1$ that is covered by chart
(1.2.2.4) but not by (1.2.2.3) is isomorphic to $\mathbb{A}^{m-3}$
and so on, until we add an affine line to $G_1$ in chart
$(1.2.2.m)$. The intersection of $G_1$ with $E_2$ is isomorphic to
$\mathbb{P}^{m-3}$. It is not so hard to see that $H(D_1\cap
E_2\cap G_1)=c_{m-2}$ (notice that the equations of (the strict
transform of) $D_1$ in chart (1.2.2.3) for instance are $x_1=0$
and $1+x_4^2+\cdots+x_m^2=0$), and from this it follows that
$H((D_1\cap G_1)^{\circ})=(w-1)c_{m-2}$ (the $w$ comes from the
$x_2$-coordinate that can be chosen freely in every chart). Now we
also have $H((E_2\cap G_1)^{\circ})=w^{m-3}+\cdots +1-c_{m-2}$ and
$H(G_1^{\circ})=w^{m-2}+\cdots
+w-(w^{m-3}+\cdots+1)-wc_{m-2}+c_{m-2}=w^{m-2}-1-(w-1)c_{m-2}$.
One gets from charts (2.1.2) to $(2.1.m)$ that $H((E_1\cap
G_1)^{\circ})=w^{m-3}+\cdots +1 -c_{m-2}$ and that $H(D_1\cap
E_1\cap G_1)=c_{m-2}$.

More conceptually, $G_1$ is a locally trivial
$\mathbb{P}^{m-3}$-bundle over the singular line on $D_1$ and
$E_1\cap G_1$ and $E_2\cap G_1$ are two fibers. Thus
$H(G_1)=(w+1)(w^{m-3}+\cdots+1)$ and $H(E_i\cap
G_1)=w^{m-3}+\cdots+1$. Furthermore, we can consider the singular
line on $D_1$ as a family of $A_1$ singularities and thus $D_1\cap
G_1$ is a family of nonsingular quadrics in $\mathbb{P}^{m-3}$.
This implies that $H(D_1\cap G_1)=(w+1)c_{m-2}$ and $H(D_1\cap
E_i\cap G_1)=c_{m-2}$. \\

In exactly the same way one finds that (for $i\in
\{2,\ldots,k-2\}$) $H(G_i^{\circ})=w^{m-2}-1-(w-1)c_{m-2}$,
$H((D_i\cap G_i)^{\circ})=(w-1)c_{m-2}$, $H((E_i\cap
G_i)^{\circ})=H((E_{i+1}\cap
G_i)^{\circ})=w^{m-3}+\cdots+1-c_{m-2}$ and $H(D_i\cap E_i\cap
G_i)=H(D_i\cap E_{i+1}\cap G_i)=c_{m-2}$.\\

\emph{Step $n+k-1$}: This step looks very much like step $n+1$. It
suffices to consider charts $(1'.2.2.1)$ to $(1'.2.2.m)$ to
compute $H(G_{k-1}^{\circ})$. One checks that $H(D_{k-1}\cap
F_1\cap G_{k-1})=H(D_{k-1}\cap F_2\cap G_{k-1})=c_{m-2}$,
$H((F_1\cap G_{k-1})^{\circ})=H((F_2\cap
G_{k-1})^{\circ})=w^{m-3}+\cdots+1-c_{m-2}$, $H((D_{k-1}\cap
G_{k-1})^{\circ})=(w-2)c_{m-2}$ and thus
$H(G_{k-1}^{\circ})=w^{m-2}+\cdots
+w-2(w^{m-3}+\cdots+1)-(w-2)c_{m-2}$. From charts $(2'.1.1.1.2)$
to $(2'.1.1.1.m)$ we get $H(D_{k-1}\cap E_{k-1} \cap
G_{k-1})=c_{m-2}$ and $H((E_{k-1}\cap
G_{k-1})^{\circ})=w^{m-3}+\cdots +1-c_{m-2}$.

A conceptual explanation like in step $n+1$ can be given here too.\\

\noindent \underline{(2) $n$ odd}\\

There are only 7 changes in comparison with the case where $n$ is
even. First remark that $ F_1\cap G_{k-1}$ and $D_{k-1}\cap
F_1\cap G_{k-1}$ are empty, but instead $H((F_1\cap
F_2)^{\circ})=c_{m-1}-c_{m-2}$ and $H(D_{k-1}\cap F_1\cap
F_2)=c_{m-2}$. The other 5 changes are the following:

\[\begin{array}{l}
H(F_1^{\circ})=b_{m-1}-b_{m-2}\\
H(F_2^{\circ})=c_m-c_{m-1}-b_{m-2}+c_{m-2}\\
H(G_{k-1}^{\circ})=w^{m-2}-1-(w-1)c_{m-2}\\
H((D_{k-1}\cap F_2)^{\circ})=b_{m-2}-c_{m-2}-1\\
H((D_{k-1}\cap G_{k-1})^{\circ})=(w-1)c_{m-2}
\end{array}\]

\noindent \textbf{3.5.} \begin{tabular}{|c|} \hline \textbf{Case}
$\mathbf{E_6}$
\\ \hline
\end{tabular}\\

\noindent We just list the results.

\setlongtables
\begin{longtable}{l}
$H(D_1^{\circ})= a_{m-2}-w-1$\\
$H(D_2^{\circ})= b_{m-1}-b_{m-2}$\\
$H(D_3^{\circ})= b_{m-1}-b_{m-2}-c_{m-1}+c_{m-2}$\\
$H(D_4^{\circ})= c_m-b_{m-2}-c_{m-1}+c_{m-2}$\\
$H(D_5^{\circ})= w^{m-2}+\cdots +w-wc_{m-2}$\\
$H((D_1\cap D_2)^{\circ})=b_{m-2}-1$\\
$H((D_1\cap D_3)^{\circ})=H((D_1\cap D_4)^{\circ})=b_{m-2}-c_{m-2}-1$\\
$H((D_1\cap D_5)^{\circ})=wc_{m-2}$\\
$H((D_2\cap D_3)^{\circ})=H((D_3\cap D_4)^{\circ})=c_{m-1}-c_{m-2}$\\
$H((D_4\cap D_5)^{\circ})=w^{m-3}+\cdots +1-c_{m-2}$\\
$H(D_1\cap D_2\cap D_3)=H(D_1\cap D_3\cap D_4)=H(D_1\cap D_4\cap
D_5)=c_{m-2}$
\end{longtable}

\newpage

\noindent \textbf{3.6.} \begin{tabular}{|c|} \hline \textbf{Cases}
$\mathbf{E_7}$ \textbf{and} $\mathbf{E_8}$
\\ \hline
\end{tabular}\\

\noindent Let us first treat the $E_8$ case. From the intersection
diagram it follows that we have to compute 47 Hodge-Deligne
polynomials (there are 12 divisors, 23 intersections of 2 divisors
and 12 intersections of 3 divisors). But there are 20 polynomials
coming from the `$D_6$ part' of the diagram that are left
unchanged here. So we will only write down the other 27.
\[\begin{array}{l}
H(B_1^{\circ})= a_{m-2}-w-1\\
H(C_1^{\circ})= a_{m-2}-b_{m-2}-w\\
H(D_1^{\circ})=H(D_2^{\circ})= a_{m-2}-2b_{m-2}+c_{m-2}-w+1\\
H(E_1^{\circ})=H(F_1^{\circ})= c_m-2b_{m-2}+c_{m-2}\\
H(H_1^{\circ})=H(I_1^{\circ})= w^{m-2}+\cdots +w-wc_{m-2}\\
H((B_1\cap C_1)^{\circ})=H((B_1\cap I_1)^{\circ})=H((C_1\cap H_1)^{\circ})=wc_{m-2}\\
H((B_1\cap D_1)^{\circ})=H((B_1\cap E_1)^{\circ})= H((C_1\cap
D_1)^{\circ})\\ \quad  =H((C_1\cap D_2)^{\circ})= H((C_1\cap
F_1)^{\circ})=H((D_1\cap D_2)^{\circ}) \\ \quad = H((D_1\cap
E_1)^{\circ})=H((D_2\cap F_1)^{\circ})=b_{m-2}-c_{m-2}-1\\
H((E_1\cap I_1)^{\circ})=H((F_1\cap H_1)^{\circ})=w^{m-3}+\cdots +1-c_{m-2}\\
H(B_1\cap C_1\cap D_1)=H(B_1\cap D_1\cap E_1)=H(B_1\cap E_1\cap
I_1)\\ \quad =H(C_1\cap D_1\cap D_2)=H(C_1\cap D_2\cap
F_1)=H(C_1\cap F_1\cap H_1)=c_{m-2}
\end{array}\]

For the $E_7$ case, we can skip all expressions involving the
divisors $B_1$ and/or $I_1$. This leaves us with 37 polynomials
and apart from the following 5, they are all the same as in the
$E_8$ case.
\[\begin{array}{l}
H(C_1^{\circ})=a_{m-2}-w-1\\
H(D_1^{\circ})=a_{m-2}-b_{m-2}-w\\
H(E_1^{\circ})=c_{m}-b_{m-2}\\
H((C_1\cap D_1)^{\circ})=H((D_1\cap E_1)^{\circ})=b_{m-2}-1 \\
\end{array}\]

$\phantom{some place}$

\section{Computation of the discrepancy coefficients}

\noindent \textbf{4.1.} In this section we compute the last data
that we need: the discrepancy coefficients. As already mentioned
in (1.4), all the two dimensional $A$-$D$-$E$'s admit a crepant
resolution, this means that all the discrepancies are 0.

For the three-dimensional case, the computations are done in
\cite{DaisRoczen}, but the authors are a bit inaccurate. Let us
again consider the case $D_n$, $n$ even, with $k=\frac{n}{2}$. The
intersection diagram is as follows:\\[-10mm]
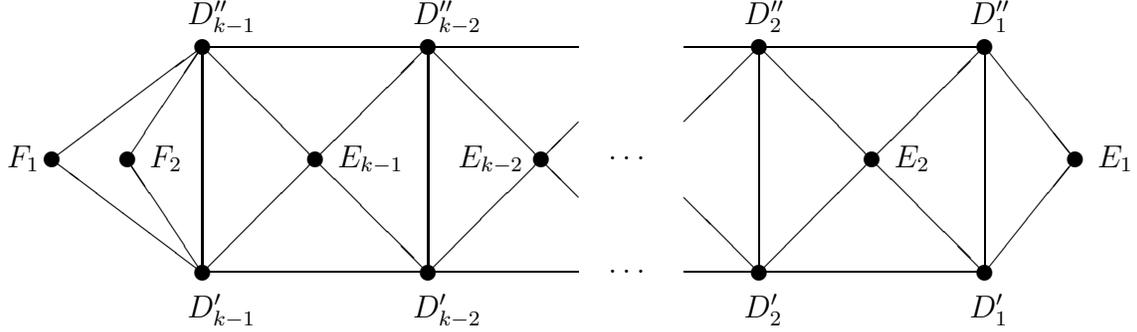
\begin{center}
\begin{picture}(160,45)
\put(27,7){\circle*{2}} \put(57,7){\circle*{2}}
\put(101,7){\circle*{2}} \put(131,7){\circle*{2}}
\put(27,37){\circle*{2}} \put(57,37){\circle*{2}}
\put(101,37){\circle*{2}} \put(131,37){\circle*{2}}
\put(42,22){\circle*{2}} \put(72,22){\circle*{2}}
\put(116,22){\circle*{2}} \put(143,22){\circle*{2}}
\put(7,22){\circle*{2}} \put(17,22){\circle*{2}}
\put(27,7){\line(1,0){50}} \put(91,7){\line(1,0){40}}
\put(27,7){\line(0,1){30}} \put(57,7){\line(0,1){30}}
\put(101,7){\line(0,1){30}} \put(131,7){\line(0,1){30}}
\put(27,7){\line(1,1){30}} \put(57,7){\line(1,1){20}}
\put(101,7){\line(1,1){30}} \put(131,7){\line(4,5){12}}
\put(101,37){\line(1,-1){30}} \put(131,37){\line(4,-5){12}}
\put(27,37){\line(1,-1){30}} \put(57,37){\line(1,-1){20}}
\put(101,7){\line(-1,1){10}} \put(101,37){\line(-1,-1){10}}
\put(27,7){\line(-2,3){10}} \put(27,7){\line(-4,3){20}}
\put(17,22){\line(2,3){10}} \put(7,22){\line(4,3){20}}
\put(27,37){\line(1,0){50}} \put(91,37){\line(1,0){40}}
\put(25,40){$D''_{k-1}$} \put(55,40){$D''_{k-2}$}
\put(99,40){$D''_2$} \put(129,40){$D''_{1}$}
\put(25,1){$D_{k-1}'$} \put(55,1){$D_{k-2}'$} \put(99,1){$D_{2}'$}
\put(129,1){$D_{1}'$} \put(1,21){$F_1$} \put(20,21){$F_2$}
\put(45,21){$E_{k-1}$} \put(61,21){$E_{k-2}$} \put(119,21){$E_2$}
\put(146,21){$E_1$} \put(81,22){$\ldots$} \put(81,7){$\ldots$}
\end{picture}
\end{center}

Compared to the higher dimensional cases, the $D_i$ fall apart
into two components $D_i'$ and $D''_{i}$, and there are no
divisors $G_i$ needed. If we denote by $\varphi:\widetilde{X}\to
X$ the log resolution, with $X$ the defining variety of the $D_n$
singularity and $\widetilde{X}$ the strict transform of $X$, then
$\varphi$ can be decomposed into $k$ birational morphisms
\[ \begin{array}{ccccccccccc} &\varphi_k & & & & & & \varphi_2& &
\varphi_1& \\
\widetilde{X}=X_k & \longrightarrow& X_{k-1}&\longrightarrow &
\cdots& \longrightarrow & X_2&\longrightarrow &X_1 &
\longrightarrow &X_0=X,
\end{array} \]
where the exceptional locus of $\varphi_1$ is $\{D'_1,D''_1\}$, of
$\varphi_i$ ($2\leq i\leq k-1$) is $\{D'_i,D''_i,E_{i-1}\}$ and of
$\varphi_k$ is $\{F_1,F_2,E_{k-1}\}$, again using the same name
for the divisors at any stage of the decomposition of $\varphi$.
We can also decompose $K_{\widetilde{X}}-\varphi^*(K_X)$ as
\[\left[\sum_{i=1}^{k-1} \varphi_k^*(\varphi_{k-1}^* \cdots
(\varphi_{i+1}^*(K_{X_i}-\varphi_i^*(K_{X_{i-1}})))\cdots)\right]
+ K_{X_k}-\varphi_k^*(K_{X_{k-1}}).\]

Dais and Roczen calculated that for instance
$\varphi_2^*(D'_1)=D'_1+D'_2+E_1$ and
$\varphi_2^*(D''_1)=D''_1+D''_2+E_1$, but $D'_1$ and $D''_1$ are
not Cartier. Their sum $D'_1+D''_1$ is Cartier and it turns out
that $\varphi_2^*(D'_1+D''_1)=D'_1+D''_1+D'_2+D''_2+E_1$ instead
of $\cdots +2E_1$. This kind of error occurs also in the following
stages for this type of singularity and also for type $D_n$, $n$
odd, and for types $E_6,E_7$ and $E_8$. In the next table, we list
the discrepancies. We use notations analogous to our notations
from section 2, but they differ from the notations in
\cite{DaisRoczen}. The coefficients
that we have corrected are in boldface.\\

\setlongtables
\begin{longtable}{l|l|c}
\multicolumn{2}{c|}{} & \\ \multicolumn{2}{c|}{Type of
singularity} &
\multicolumn{1}{c}{Discrepancy} \\
\multicolumn{2}{c|}{} & \\ \hline & & \\ $A_n$ & $\begin{array}{l} n \text{ even } \\
n=2k \\ k\geq 1  \end{array}$
& $\displaystyle{\sum_{i=1}^k i D_i  + (n+2)D_{k+1}}$ \\
& & \\ \cline{2-3} \newpage \cline{2-3} & &  \\ & $\begin{array}{l} n \text{ odd } \\
n=2k-1 \\ k\geq 1  \end{array}$ &
$\displaystyle{\sum_{i=1}^k i D_i} $ \\ & &  \\ \hline & &  \\
$D_n$ & $\begin{array}{l} n \text{ even } \\
n=2k \\ k\geq 2  \end{array}$ & $\displaystyle{\sum_{i=1}^{k-1}
\bigl(iD'_i+iD''_i+\mathbf{2i}E_i \bigr)+ \mathbf{k}F_1+\mathbf{k}F_2}$  \\ & &  \\
\cline{2-3} & &  \\ & $\begin{array}{l} n \text{ odd } \\
n=2k+1 \\ k\geq 2  \end{array}$ & $\displaystyle{\sum_{i=1}^{k-1}
\bigl(iD'_i+iD''_i+\mathbf{2i}E_i
\bigr)+ \mathbf{k}F_1+\mathbf{2k}F_2}$ \\ & &  \\ \hline \multicolumn{2}{c|}{} &  \\
\multicolumn{2}{c|}{$E_6$}
& $D'_1+D''_1+\mathbf{2}D_2+\mathbf{4}D_3+\mathbf{6}D_4 $  \\
\multicolumn{2}{c|}{} &   \\ \hline \multicolumn{2}{c|}{} &    \\
\multicolumn{2}{c|}{$E_7$} &
$C'_1+C''_1+2D'_1+2D''_1+4D'_2+4D''_2$
\\ \multicolumn{2}{c|}{} & $+\mathbf{3}E_1+\mathbf{7}E_2+\mathbf{6}F_1+
\mathbf{5}F_2$ \\ \multicolumn{2}{c|}{} &    \\ \hline \multicolumn{2}{c|}{} & \\
\multicolumn{2}{c|}{$E_8$} &
$B'_1+B''_1+2C'_1+2C''_1+4D'_1+4D''_1+7D'_2+7D''_2$\\
\multicolumn{2}{c|}{} & $+\mathbf{6}E_1+\mathbf{12}E_2
+\mathbf{10}F_1+ \mathbf{8}F_2$ \\
\multicolumn{2}{c|}{} & \\
\end{longtable}

\noindent \textbf{Remark.} Dais and Roczen used their results to
contradict a conjecture of Batyrev about the range of the
string-theoretic index (see \cite[Conjecture 5.9]{Batyrev},
\cite[Remark 1.9]{DaisRoczen}). Luckily, this follows already from
the formulae for the $A$ case, to which we do not correct
anything. We will only simplify their formulae in this case.\vspace{3.2mm}

\noindent \textbf{4.2.} Now we consider the higher dimensional
case. As an example, we will calculate the discrepancy coefficient
of the divisor $E_i$ for an $(m-1)$-dimensional $D_n$ singularity,
where $n$ is even, $i\in \{1,\ldots,k-1\}$ and $m\geq 5$. Let $X$
be the defining variety $\{x_1^{n-1}+x_1x_2^2+x_3^2+\cdots + x_m^2
= 0\}\subset\mathbb{A}^m$, and let $\varphi: \widetilde{X}\to X$
be the log resolution constructed in section 2. We take a
coordinate chart that covers a piece of $E_i$; in the notation of
section 2, this could be for example
chart\begin{tabular}[t]{c}($\underbrace{1.1\ldots 1}$.2.3)\\ $i-1$
times
\end{tabular}describing an open set $U\subset \widetilde{X}$: \[
y_1^{2k-2i+1}y_2^{2k-2i-1}y_3^{4k-4i-2}+y_1y_2+1+y_4^2+\cdots
+y_m^2=0.\] In this chart, $y_1=0$ gives a local equation for
divisor $D_{i-1}$, $y_2=0$ for $D_i$ and $y_3=0$ for our divisor
$E_i$. The map $\varphi:U\to X$ can be found from the resolution
process. Here it will be
\[\varphi(y_1,\ldots,y_m)=(y_1y_2y_3^2,y_1^{i-1}y_2^{i}y_3^{2i-1},
y_1^{i-1}y_2^{i}y_3^{2i},y_1^{i-1}y_2^{i}y_3^{2i}y_4, \ldots,
y_1^{i-1}y_2^{i}y_3^{2i}y_m).\]

The section $\frac{dx_1\wedge \ldots \wedge dx_{m-1}}{2x_m}$ is
locally a generator of the sheaf $\mathcal{O}_X(K_X)$
($2x_m=\frac{\partial\, f}{\partial\, x_m}$, where $f$ is the
equation of $X$) and we have to compare its pull-back under
$\varphi$ with the generator $\frac{dy_1\wedge \ldots \wedge
dy_{m-1}}{2y_m}$ of
$\mathcal{O}_{\widetilde{X}}(K_{\widetilde{X}})|_U$. We have \[
\varphi^*(\frac{dx_1\wedge \ldots \wedge
dx_{m-1}}{2x_m})=y_1^{(i-1)(m-3)}y_2^{i(m-3)}y_3^{2i(m-3)}
\frac{dy_1\wedge \ldots \wedge dy_{m-1}}{2y_m},\] which learns us
that the discrepancy coefficient of $E_i$ is $2i(m-3)$. And we get
the discrepancy coefficient of $D_i$ for free, it is $i(m-3)$. In
general, the following can be proven by this kind of calculations.\\

\noindent \textbf{Proposition.} \textsl{For all divisors that are created
after a point blow-up, except for divisor $D_{\frac{n}{2}+1}$ in
the $A_n$ ($n$ even) case, the discrepancy coefficient is ($m-3$)
times the coefficient of the corresponding divisor(s) in the
three-dimensional case (see the table in (4.1)).}\\

What about the other divisors$\,$? They are all created after
blowing up a nonsingular surrounding variety in a point (case
$A_n$, $n$ even) or a line (other cases). We consider again the
case of a $D_n$ singularity, with $n$ even. Denote by $X^{(i)}$
the variety obtained after $n+i$ steps in the resolution process
of section 2 ($i\in \{0,\ldots, k-2\}$). The log resolution
$\varphi:\widetilde{X}\to X $ can be decomposed as follows:
\[ \begin{array}{ccccccc} &\chi^{(i+1)} &  & \varphi^{(i+1)}& &
\psi^{(i)}& \\
\widetilde{X} & \longrightarrow& X^{(i+1)}&\longrightarrow &
X^{(i)}&\longrightarrow &X,
\end{array} \]
where $\varphi^{(i+1)}$ is the blow-up of the singular line on the
divisor $D_{i+1}\subset X^{(i)}$ and where $\chi^{(i+1)}$ and
$\psi^{(i)}$ are compositions of other blow-ups. Notice that all
the singular lines on $X^{(0)}$ are disjoint. Thus, to compute the
discrepancy coefficient of $G_{i+1}$, it suffices to look at its
coefficient in
$K_{X^{(i+1)}}-(\psi^{(i)}\circ\varphi^{(i+1)})^*(K_X)$. This is
equal to
\[ K_{X^{(i+1)}}-
(\varphi^{(i+1)})^*((\psi^{(i)})^*(K_X)-K_{X^{(i)}})-
(\varphi^{(i+1)})^*(K_{X^{(i)}}).
\]

It follows from \cite[p.608]{GriffithsHarris} that the last term
is $-K_{X^{(i+1)}}+(m-3)G_{i+1}$ ($X^{(i)}$ is nonsingular$\,$!).
And in the second term we only get a nonzero coefficient for
$G_{i+1}$ from $-(\varphi^{(i+1)})^*(-(i+1)(m-3)D_{i+1})$ (this
follows from \cite[p.605]{GriffithsHarris}, and the exact
coefficient is $2(i+1)(m-3)$ because the multiplicity of a generic
point of the singular line on $D_{i+1}$ is 2). This gives us
$2(i+1)(m-3)+(m-3)=(2i+3)(m-3)$ as discrepancy coefficient for
$G_{i+1}$. In all other cases where we blow up in a line, the
multiplicity of a generic point of the singular line will also be
2 and thus we have the following proposition.\\

\noindent \textbf{Proposition.} \textsl{For all divisors that are created
after a blow-up in a singular line of another divisor $D$, the
discrepancy coefficient is \[2(\text{discrepancy coefficient of
}D) + (m-3).\]}\indent The reader may check that the same arguments give $(n+1)(m-3)+1$
as coefficient for $D_{\frac{n}{2}+1}$ in the case $A_n$, $n$
even.\\

\section{Formulae for the contribution of an $A$-$D$-$E$
singularity to the stringy $E$-function and application to
Batyrev's conjecture}

\noindent \textbf{5.1.} Let $X$ be a defining variety of an
$A$-$D$-$E$ singularity; hence $X$ is a hypersurface in
$\mathbb{A}^m$ ($m\geq 3$) with a singular point in the origin. By
the contribution of the singular point to the stringy
$E$-function, we mean $E_{st}(X)-H(X\setminus\{0\})$ (see (3.2)).
Before stating the formulae, we first remark that we have to make
a distinction between $m$ even and $m$ odd, because the required
Hodge-Deligne polynomials depend on the parity of the dimension.\\

\noindent \textbf{Theorem.} \textsl{The contributions of the
($m-1$)-dimensional $A$-$D$-$E$ singularities ($m\geq 3$) are
given in the following tables (where sums like $\sum_{i=2}^k$ must
be interpreted as $0$ for $k=1$).}
\setlongtables
\begin{longtable}{l|l|c}
\multicolumn{2}{c|}{} & \\ \multicolumn{2}{c|}{Type of
singularity} &
\multicolumn{1}{c}{Contribution of singular point for odd $m$} \\
\multicolumn{2}{c|}{} & \\ \hline & & \\  $A_n$ &
$\begin{array}{l} n \text{ even } \\ n=2k \\ k\geq 1  \end{array}$
& $\displaystyle{1+\frac{(w-1)}{(w^{(2k+1)(m-3)+2}-1)}\left(
 \sum_{i=2}^{k+1} w^{(k+i)(m-3)+2} \right.}$\\
& & $\displaystyle{\left. +\sum_{i=1}^k
w^{(k+i)(m-3)+\frac{m+1}{2}}+\sum_{i=1}^k
w^{i(m-3)+\frac{m-1}{2}}+ \sum_{i=1}^k w^{i(m-3)+1}
 \right)}$ \\
& & \\ \cline{2-3} \newpage  \cline{2-3} & & \\ & $\begin{array}{l} n \text{ odd } \\
n=2k-1 \\ k\geq 1  \end{array}$ &
$\displaystyle{1+\frac{(w-1)}{(w^{k(m-3)+1}-1)}\left(
\sum_{i=1}^{k} w^{i(m-3)+1}+ \sum_{i=1}^{k-1}
w^{i(m-3)+\frac{m-1}{2}}
 \right)} $ \\ & &  \\ \hline  & &  \\
$D_n$ & $\begin{array}{l} n \text{ even } \\ n=2k \\ k\geq 2
\end{array}$ & $\displaystyle{1+\frac{(w-1)}{(w^{(2k-1)(m-3)+1}-1)}\left(
\sum_{i=1}^{2k-1} w^{i(m-3)+1}+ w^{k(m-3)+1}
 \right)}$  \\ & &  \\
\cline{2-3} & &  \\ & $\begin{array}{l} n \text{ odd } \\
n=2k+1 \\ k\geq 2  \end{array}$ &
$\displaystyle{1+\frac{(w-1)}{(w^{2k(m-3)+1}-1)}\biggl(
 \sum_{i=1}^{2k} w^{i(m-3)+1} +
w^{k(m-3)+\frac{m-1}{2}}
 \biggr)}$ \\ & &  \\ \hline \multicolumn{2}{c|}{} &  \\
\multicolumn{2}{c|}{$E_6$} &
$\displaystyle{1+\frac{(w-1)}{(w^{6m-17}-1)}\Bigl(w^{6m-17}+w^{4m-11}+w^{3m-8} \Bigr.}$\\
\multicolumn{2}{c|}{} &
$\displaystyle{\Bigl.+w^{m-2}+w^{\frac{9m-25}{2}}+w^{\frac{5m-13}{2}}
  \Bigr)  }$\\
\multicolumn{2}{c|}{} &   \\ \hline \multicolumn{2}{c|}{} &    \\
\multicolumn{2}{c|}{$E_7$} &
$\displaystyle{1+\frac{(w-1)}{(w^{9m-26}-1)}\left(w^{9m-26}+
w^{7m-20}+w^{6m-17}+w^{5m-14}\right.}$\\
\multicolumn{2}{c|}{} & $\displaystyle{\left.+w^{4m-11}
+w^{3m-8}+w^{m-2} \right)}$
\\ \multicolumn{2}{c|}{} &
 \\ \hline \multicolumn{2}{c|}{} & \\
\multicolumn{2}{c|}{$E_8$} &
$\displaystyle{1+\frac{(w-1)}{(w^{15m-44}-1)}\left(w^{15m-44}+w^{12m-35}+w^{10m-29}+
w^{9m-26}\right.}$\\
\multicolumn{2}{c|}{} & $\displaystyle{\left.+
w^{7m-20}+w^{6m-17}+w^{4m-11}+w^{m-2} \right)}$
 \\
\multicolumn{2}{c|}{} & \\
\end{longtable}

$\phantom{some place}$

\setlongtables
\begin{longtable}{l|l|c}
\multicolumn{2}{c|}{} & \\ \multicolumn{2}{c|}{Type of
singularity} &
\multicolumn{1}{c}{Contribution of singular point for even $m$} \\
\multicolumn{2}{c|}{} & \\ \hline & & \\  $A_n$ &
$\begin{array}{l} n \text{ even } \\ n=2k \\ k\geq 1  \end{array}$
& $\displaystyle{1+\frac{(w-1)}{(w^{(2k+1)(m-3)+2}-1)}\left(
 \sum_{i=2}^{k+1} w^{(k+i)(m-3)+2} + \sum_{i=1}^k w^{i(m-3)+1}
 \right)}$ \\
& & \\ \cline{2-3} \newpage \cline{2-3} & &  \\ & $\begin{array}{l} n \text{ odd } \\
n=2k-1 \\ k\geq 1  \end{array}$ &
$\displaystyle{1+\frac{(w-1)}{(w^{k(m-3)+1}-1)}\left(
\sum_{i=1}^{k} w^{i(m-3)+1}+ w^{\frac{m}{2}-1}
 \right)} $ \\ & &  \\ \hline  & &  \\
$D_n$ & $\begin{array}{l} n \text{ even } \\ n=2k \\ k\geq 2
\end{array}$ & $\displaystyle{1+\frac{(w-1)}{(w^{(2k-1)(m-3)+1}-1)}\left(
\sum_{i=1}^{2k-1} w^{i(m-3)+1}+ w^{k(m-3)+1}\right.}$\\ & &
$\displaystyle{\left.+\sum_{i=0}^{k-2} w^{(k+i)(m-3)+\frac{m}{2}}
+ \sum_{i=0}^{k-1} w^{i(m-3)+\frac{m}{2}-1}+w^{\frac{m}{2}-1}
 \right)}$  \\ & &  \\
\cline{2-3} & &  \\ & $\begin{array}{l} n \text{ odd } \\
n=2k+1 \\ k\geq 2  \end{array}$ &
$\displaystyle{1+\frac{(w-1)}{(w^{2k(m-3)+1}-1)}\biggl(
 \sum_{i=1}^{2k} w^{i(m-3)+1}+ \sum_{i=1}^{k-1}
w^{(k+i)(m-3)+\frac{m}{2}} \biggr.}$\\ & & $\displaystyle{\biggl.
+\sum_{i=0}^{k-1} w^{i(m-3)+\frac{m}{2}-1}
 \biggr)}$ \\ & &  \\ \hline \multicolumn{2}{c|}{} &  \\
\multicolumn{2}{c|}{$E_6$} &
$\displaystyle{1+\frac{(w-1)}{(w^{6m-17}-1)}\Bigl(w^{6m-17}+w^{4m-11}+w^{3m-8} \Bigr.}$\\
\multicolumn{2}{c|}{} &
$\displaystyle{\Bigl.+w^{m-2}+w^{\frac{11m-30}{2}}+w^{\frac{3m-8}{2}}
  \Bigr)  }$\\
\multicolumn{2}{c|}{} &   \\ \hline \multicolumn{2}{c|}{} &    \\
\multicolumn{2}{c|}{$E_7$} &
$\displaystyle{1+\frac{(w-1)}{(w^{9m-26}-1)}\Bigl(w^{9m-26}+
w^{7m-20}+w^{6m-17}+w^{5m-14}\Bigr.}$\\
\multicolumn{2}{c|}{} & $\displaystyle{\Bigl.
+w^{4m-11}+w^{3m-8}+w^{m-2}
+w^{\frac{17m-48}{2}}+w^{\frac{15m-42}{2}}+w^{\frac{11m-30}{2}}\Bigr.}$\\
\multicolumn{2}{c|}{} & $\displaystyle{\Bigl.+w^{\frac{9m-26}{2}}
+w^{\frac{5m-14}{2}}+w^{\frac{3m-8}{2}}+w^{\frac{m-2}{2}}\Bigr)}$
\\ \multicolumn{2}{c|}{} &
 \\ \hline \multicolumn{2}{c|}{} & \\
\multicolumn{2}{c|}{$E_8$} &
$\displaystyle{1+\frac{(w-1)}{(w^{15m-44}-1)}\Bigl(w^{15m-44}+w^{12m-35}+w^{10m-29}+
w^{9m-26}\Bigr.}$\\
\multicolumn{2}{c|}{} & $\displaystyle{\Bigl.+
w^{7m-20}+w^{6m-17}+w^{4m-11}+w^{m-2}+w^{\frac{29m-84}{2}}+w^{\frac{27m-78}{2}}
\Bigr.}$\\
\multicolumn{2}{c|}{} & $\displaystyle{\Bigl.
+w^{\frac{23m-66}{2}}+w^{\frac{17m-48}{2}}+
w^{\frac{15m-44}{2}}+w^{\frac{9m-26}{2}}+w^{\frac{5m-14}{2}}+w^{\frac{3m-8}{2}}
 \Bigr)}$
 \\
\multicolumn{2}{c|}{} & \\
\end{longtable}

\newpage

\noindent \textbf{Proof:} 
\begin{itemize}
\item Let us first consider the case where $m\geq 5$. We will
focus again on the singularity of type $D_n$ for $n=2k$ and also
for even $m$. All the other cases are completely analogous. We
just insert the data from sections 2, 3 and 4 in the defining
formula of the stringy $E$-function and we find the following
formula for the contribution of the singularity:
\begin{equation*}
\scriptsize{
\begin{split}
&
\frac{(w^{m-1}-w^2+w^{\frac{m+2}{2}}-w^{\frac{m}{2}})}{(w^{m-2}-1)}
+ \sum_{i=2}^{k-1}
\frac{(w^{m-2}-w+w^{\frac{m}{2}}-w^{\frac{m-2}{2}}
)(w-1)}{(w^{i(m-3)+1}-1)} + \frac{(w^{m-2})(w-1)}{(w^{2m-5}-1)}
\\ & \  + \sum_{i=2}^{k-1}
\frac{(w^{m-2}-w^{m-3}-w^{\frac{m-2}{2}}+w^{\frac{m-4}{2}}
)(w-1)}{(w^{2i(m-3)+1}-1)} +
\frac{2w^{m-2}(w-1)}{(w^{k(m-3)+1}-1)}\\ & \  + \sum_{i=1}^{k-2}
\frac{(w^{m-2}-w^{m-3}-w^{\frac{m-2}{2}}+
w^{\frac{m-4}{2}})(w-1)}{(w^{(2i+1)(m-3)+1}-1)} +
\frac{(w^{m-2}-2w^{m-3}-w^{\frac{m-2}{2}}+2w^{\frac{m-4}{2}})
(w-1)}{(w^{(2k-1)(m-3)+1}-1)}\\ & \  + \sum_{i=1}^{k-2}
\frac{(w^{m-2}-w+w^{\frac{m}{2}}-w^{\frac{m-2}{2}}
)(w-1)}{(w^{i(m-3)+1}-1)(w^{(i+1)(m-3)+1}-1)} +
\frac{(w^{m-2}-w+w^{\frac{m}{2}}-w^{\frac{m-2}{2}}
)(w-1)}{(w^{m-2}-1)(w^{2m-5}-1)}\\ & \    + \sum_{i=2}^{k-1}
\frac{(w^{m-3}-1+w^{\frac{m-2}{2}}-w^{\frac{m-4}{2}})
(w-1)^2}{(w^{i(m-3)+1}-1)(w^{2i(m-3)+1}-1)}+ \sum_{i=1}^{k-2}
\frac{(w^{m-3}-1+w^{\frac{m-2}{2}}-w^{\frac{m-4}{2}})
(w-1)^2}{(w^{i(m-3)+1}-1)(w^{(2i+2)(m-3)+1}-1)}\\ &\  +
\frac{2(w^{m-2}-w+w^{\frac{m}{2}}-w^{\frac{m-2}{2}}
)(w-1)}{(w^{(k-1)(m-3)+1}-1)(w^{k(m-3)+1}-1)} + \sum_{i=1}^{k-2}
\frac{(w^{m-3}-1+w^{\frac{m-2}{2}}-w^{\frac{m-4}{2}})
(w-1)^2}{(w^{i(m-3)+1}-1)(w^{(2i+1)(m-3)+1}-1)}\\ & \  +
\frac{(w^{m-2}-2w^{m-3}-w+2+w^{\frac{m}{2}}-3w^{\frac{m-2}{2}}+2w^{\frac{m-4}{2}})
(w-1)}{(w^{(k-1)(m-3)+1}-1)(w^{(2k-1)(m-3)+1}-1)}\\ & \  +
\sum_{i=1}^{k-1} \frac{(w^{m-3}-w^{\frac{m-4}{2}})
(w-1)^2}{(w^{2i(m-3)+1}-1)(w^{(2i+1)(m-3)+1}-1)} +
\sum_{i=1}^{k-2} \frac{(w^{m-3}-w^{\frac{m-4}{2}})
(w-1)^2}{(w^{(2i+2)(m-3)+1}-1)(w^{(2i+1)(m-3)+1}-1)}\\ & \  +
\frac{2(w^{m-3}-w^{\frac{m-4}{2}})
(w-1)^2}{(w^{k(m-3)+1}-1)(w^{(2k-1)(m-3)+1}-1)}\\& \  +
\sum_{i=1}^{k-2}
\frac{(w^{m-3}-1+w^{\frac{m-2}{2}}-w^{\frac{m-4}{2}})(w-1)^2}{(w^{i(m-3)+1}-1)
(w^{(i+1)(m-3)+1}-1)(w^{(2i+2)(m-3)+1}-1)} \\ & \  +
\sum_{i=1}^{k-1}
\frac{(w^{m-3}-1+w^{\frac{m-2}{2}}-w^{\frac{m-4}{2}})(w-1)^2}{(w^{i(m-3)+1}-1)
(w^{2i(m-3)+1}-1)(w^{(2i+1)(m-3)+1}-1)} \\ &\  + \sum_{i=1}^{k-2}
\frac{(w^{m-3}-1+w^{\frac{m-2}{2}}-w^{\frac{m-4}{2}})(w-1)^2}{(w^{i(m-3)+1}-1)
(w^{(2i+2)(m-3)+1}-1)(w^{(2i+1)(m-3)+1}-1)} \\ &\  +
\frac{2(w^{m-3}-1+w^{\frac{m-2}{2}}-w^{\frac{m-4}{2}})(w-1)^2}{
(w^{(k-1)(m-3)+1}-1) (w^{k(m-3)+1}-1)(w^{(2k-1)(m-3)+1}-1)}.
\end{split}}
\end{equation*}

The terms correspond to the following pieces of the exceptional
locus (in that order):
\begin{equation*}
\begin{split}
&
D_1^{\circ},D_i^{\circ},E_1^{\circ},E_i^{\circ},F_i^{\circ},G_i^{\circ},
G_{k-1}^{\circ}, (D_i\cap D_{i+1})^{\circ},(D_1\cap E_1)^{\circ},
(D_i\cap E_i)^{\circ},\\ & \  (D_i\cap
E_{i+1})^{\circ},(D_{k-1}\cap F_i)^{\circ},(D_i\cap
G_i)^{\circ},(D_{k-1}\cap G_{k-1})^{\circ},(E_i\cap G_i)^{\circ},
\\ & \  (E_{i+1}\cap G_i)^{\circ},(F_i\cap G_{k-1})^{\circ},D_i\cap
D_{i+1}\cap E_{i+1},D_i\cap E_i\cap G_i,\\ & \  D_i\cap
E_{i+1}\cap G_i,D_{k-1}\cap F_i\cap G_{k-1}.
\end{split}
\end{equation*}

By a very long but easy calculation, it can be proved by induction
on $k$ that we indeed get the requested formula. We remark here
that we have done the computations for $m\geq 5$, for $m=4$ and
for $m=3$ separately, and then noticed that the formulae for
$m\geq 5$ are correct in the other cases too. \item We can now
explain why these formulae are also valid for $m=4$. For the $A_n$
case, this is not a surprise, since the intersection diagram for
$m=4$ is the same as for $m\geq 5$.

For the other cases, consider for example a singularity of type
$D_n$, $n$ even. The blow-ups in the singular lines on the
divisors $D_i$ in the higher dimensional case correspond here to
blow-ups in the intersections $D'_i \cap D''_i$. Performing these
unnecessary extra blow-ups yields just another log resolution, and
the formula for the contribution of the singularity for that log
resolution will be exactly the evaluation of the formula from the
first part of the proof for $m=4$ (notice for instance that the
Hodge-Deligne polynomial for $D_i^{\circ}$ becomes $2w^2-2w$ for
$m=4$ and the Hodge-Deligne polynomials for $(D'_i)^{\circ}$ and
$(D''_i)^{\circ}$ will both be $w^2-w$). \item For $m=3$ it can be
checked easily that the formulae are correct but again we give a
more conceptual explanation. Compared with the higher dimensional
case, all divisors except the last one split into two (distinct)
components in the $A_n$ case, for odd $n$. This is consistent with
the Hodge-Deligne polynomials from (3.3), evaluated for $m=3$. For
even $n$, we must notice that the last blow-up is unnecessary for
surfaces; performing it anyway does not yield a crepant resolution
any more (the last divisor has discrepancy coefficient 1, as it
should be, according to (4.2)). This last divisor is irreducible
and the first $\frac{n}{2}$ blow-ups each add two components to
the exceptional locus (compare this with (3.3) again).

For the $D_n$ case, the analogue of blowing up in a singular line
on a divisor $D_i$ would be to blow up in $D_i$ itself, because it
is just a line for $m=3$. Such a blow-up is an isomorphism, and
the result is that the divisors $D_i$ are renamed as $G_i$. As
intersection diagram one finds the same as in the higher
dimensional case, but without the divisors $D_i$. To be able to
compare this to (3.4), we must notice that it is logical to set
$a_1=w+1$, $c_1=0$ and $b_1=1$ in (3.1). Then indeed all
Hodge-Deligne polynomials that describe a piece of a divisor $D_i$
are 0 in (3.4) for $m=3$. For the $E$ cases the same sort of
arguments apply.\hfill $\blacksquare$
\end{itemize}

\noindent \textbf{5.2.} From now on, let $X$ be a projective
algebraic variety with at most (a finite number of) $A$-$D$-$E$
singularities. Since the next results are trivial for surfaces, we
will assume that $\dim X\geq 3$.\\

\noindent \textbf{Proposition.} \textsl{The stringy $E$-function of $X$ is
a polynomial if and only if $\dim X = 3$ and $X$ has singularities
of type $A_n$ ($n$ odd) and/or $D_n$ ($n$ even). }\\

\noindent \textbf{Proof:} It follows from theorem (5.1) that the contributions of the singular points for
$m\geq 5$ can be written in the following form:
\[ 1 + \frac{w^2(w^{\alpha}+a_{\alpha-1}w^{\alpha-1}+\cdots +
a_0)}{w^{\alpha+1}+w^{\alpha}+\cdots +1 },\] where $\alpha \in
\mathbb{Z}_{>0}$ and all $a_i \in \mathbb{Z}_{\geq 0}$. Such
expressions or finite sums of such expressions can never be
polynomials. For $m=4$ the contributions are given in the
following table. \setlongtables
\begin{longtable}{l|l|c}
\multicolumn{2}{c|}{$_{\displaystyle{\text{Type of
singularity}}}$} &
\multicolumn{1}{c}{$_{\displaystyle{\text{Contribution of singular point}}}$} \\
\multicolumn{2}{c|}{} & \\ \hline & & \\  $A_n$ &
$\begin{array}{l} n \text{ even } \\ n=2k   \end{array}$
& $1+\frac{w^2(w^{2k+2}-w^{k+2}+w^k-1)}{w^{2k+3}-1}$ \\
& & \\ \cline{2-3} & &  \\ & $\begin{array}{l} n \text{ odd } \\
n=2k-1   \end{array}$ &
$w+1 $ \\ & &  \\ \hline & &  \\
$D_n$ & $\begin{array}{l} n \text{ even } \\ n=2k
\end{array}$ & $2w+1$  \\ & &  \\
\cline{2-3} & &  \\ & $\begin{array}{l} n \text{ odd } \\
n=2k+1   \end{array}$ &
$w+1+\frac{w^2(w^{2k}-w^{k+1}+w^{k-1}-1)}{w^{2k+1}-1}$ \\ & &  \\
 \hline \multicolumn{2}{c|}{} &  \\
\multicolumn{2}{c|}{$E_6$} &
$1+\frac{w^2(2w^6-2w^5+w^4-w^2+2w-2)}{w^7-1}$\\
\multicolumn{2}{c|}{} &   \\ \hline \multicolumn{2}{c|}{} &    \\
\multicolumn{2}{c|}{$E_7$} & $w+1+\frac{w^2(w^4-w^3+w-1)}{w^5-1}$
\\ \multicolumn{2}{c|}{} &
 \\  \hline \multicolumn{2}{c|}{} & \\
\multicolumn{2}{c|}{$E_8$} &
$1+\frac{w^2(2w^7-w^6-w^5+2w^4-2w^3+w^2+w-2)}{w^8-1}$
 \\
\multicolumn{2}{c|}{} & \\
\end{longtable}
There are exactly two contributions that are polynomials and one
sees again that adding a finite number of the non-polynomial
expressions never gives a polynomial.\hfill $\blacksquare$\\

\newpage

\noindent \textbf{Theorem.} \textsl{Let $X$ be a three-dimensional
projective variety with at most singularities of type $A_n$ ($n$
odd) and/or $D_n$ ($n$ even). Then the stringy Hodge numbers of
$X$ are nonnegative.}\\

\noindent \textbf{Proof:} Let us first consider the case where $X$ has one singularity
of type $A_n$ ($n$ odd). Denote by $X_{ns}$ the nonsingular part
of $X$, and let $\varphi:\widetilde{X}\to X$ be the log resolution
as constructed in section 2. Then the stringy $E$-function of $X$
will be $E_{st}(X)=H(X_{ns})+uv+1$ and the Hodge-Deligne
polynomial of $\widetilde{X}$ is
$H(\widetilde{X})=H(X_{ns})+\frac{n+1}{2}(uv)^2+\frac{n+3}{2}(uv)+1$.
The exceptional locus counts $\frac{n+1}{2}$ components
$D_1,\ldots,D_{\frac{n+1}{2}}$ whose classes in
$H^2(\widetilde{X},\mathbb{C})$ are linearly independent. This can
be seen as follows. We embed $\widetilde{X}$ in a $\mathbb{P}^N$
and we intersect with a suitable hyperplane $Y$. Thanks to
Grauert's contractibility criterion the intersection matrix of the
curves $D_1\cap Y,\ldots,D_{\frac{n+1}{2}}\cap Y$ is negative
definite, and thus the classes of these curves are linearly
independent in $H^2(\widetilde{X}\cap Y,\mathbb{C})$. The weak
Lefschetz theorem implies then that the classes of
$D_1,\ldots,D_{\frac{n+1}{2}}$ are linearly independent in
$H^2(\widetilde{X},\mathbb{C})$. Actually, these classes are all
contained in $H^{1,1}(\widetilde{X})$
(\cite[p.163]{GriffithsHarris}). This means that
$h^{1,1}(\widetilde{X})=h^{2,2}(\widetilde{X})\geq \frac{n+1}{2}$.
Thus the coefficients of $(uv)^2$ and $uv$ in $H(X_{ns})$ are
$\geq 0$ and $\geq -1$ respectively. Note also that the constant
term of $H(X_{ns})$ will be zero and for all other coefficients
$a_{p,q}$ of $u^pv^q$ in $H(X_{ns})$, $(-1)^{p+q}a_{p,q}$ will be
$\geq 0$, since this is the case in $H(\widetilde{X})$. This
implies that the stringy Hodge numbers of $X$ are nonnegative.

If $X$ has one singularity of type $D_n$ ($n$ even), we can choose
to start from the log resolution constructed by Dais and Roczen
(\cite[Section 2]{DaisRoczen}) which yields
$H(\widetilde{X})=H(X_{ns})+\frac{3n-2}{2}(uv)^2+\frac{3n+2}{2}(uv)+1$
with $\frac{3n-2}{2}$ components in the exceptional locus or we
can use the log resolution analogous to section 2, which gives
$H(\widetilde{X})=H(X_{ns})+(2n-2)(uv)^2+2n(uv)+1$ with $2n-2$
components in the exceptional locus and then apply the same
argument.

It is clear that nothing essential changes when there is more than
one singularity. \hfill $\blacksquare$\\

\noindent \textbf{Example.} Consider the variety $X=\{xyz+t^3+w^3=0\} \subset \mathbb{P}^4$,
where we use coordinates $(x,y,z,t,w)$. It is clear that the
points $(1,0,0,0,0),(0,1,0,0,0)$ and $(0,0,1,0,0)$ are
three-dimensional $D_4$ singularities. Thus, their contribution to
the stringy $E$-function of $X$ is $3(2w+1)$. To calculate the
Hodge-Deligne polynomial of $X$, we divide $X$ in three locally
closed pieces:
\[ X = (X\cap \{x\neq 0, y\neq 0\})\sqcup (X\cap \{x\neq 0,y=0\})
\sqcup (X\cap \{x=0\}).\] The Hodge-Deligne polynomial of the
first piece is just $(w-1)w^2$ since $y,z,t,w$ have become affine
coordinates and $y,t,w$ can be chosen freely, with $y\neq 0$. The
second piece consists of three planes in $\mathbb{A}^3$,
intersecting in a line and has Hodge-Deligne polynomial
$3(w^2-w)+w$ and the third piece are three planes in
$\mathbb{P}^3$, intersecting in a line, with contribution
$3w^2+w+1$. Thus $H(X)=w^3+5w^2-w+1$ and $H(X_{ns})=w^3+5w^2-w-2$.
It follows that the stringy $E$-function of $X$ is equal to
$w^3+5w^2+5w+1$ and that the stringy Hodge numbers of $X$ are
nonnegative.\\

\noindent \textbf{Acknowledgement:} I wish to thank Professor J.
Steenbrink for his contribution to the proof of Theorem 5.2.\\

\footnotesize{

\end{document}
\begin{thebibliography}{DL3}

\bibitem[Ba1]{Batyrev} V. Batyrev, \emph{Stringy Hodge numbers
of varieties with Gorenstein canonical singularities}, Proc.
Taniguchi Symposium 1997, In `Integrable Systems and Algebraic
Geometry, Kobe/Kyoto 1997', World Sci. Publ. (1999), 1-32.

\bibitem[Ba2]{Batyrev2} V. Batyrev, \emph{Non-Archimedian integrals and
stringy Euler numbers of log terminal pairs}, J. Europ. Math. Soc.
\textbf{1} (1999), 5-33.

\bibitem[Da]{Dais} D. Dais, \emph{On the string-theoretic Euler number
of a class of absolutely isolated singularities}, Manuscripta
Math. \textbf{105} (2001), 143-174.

\bibitem[DR]{DaisRoczen} D. Dais and M. Roczen, \emph{On the
string-theoretic Euler number of 3-dimensional A-D-E
singularities}, Adv. Geom. \textbf{1} (2001), 373-426.

\bibitem[GH]{GriffithsHarris} P. Griffiths and J. Harris, \emph{Principles
of algebraic geometry}, John Wiley and Sons (1978).

\bibitem[Re]{Reid} M. Reid, \emph{Young person's guide to
canonical singularities}, Algebraic Geometry Bowdoin 1985, Proc.
Sympos. Pure Math., vol. 46 (1987), 345-414.

\bibitem[Sr]{Srinivas} V. Srinivas, \emph{The Hodge Characteristic},
preprint.

\bibitem[Ve1]{Veys} W. Veys, \emph{Arc spaces, motivic integration
and stringy invariants}, to appear in Adv. Stud. Pure Math., Proceedings of "Singularity Theory and its applications, Sapporo, 16-25 september 2003".

\bibitem[Ve2]{Veys2} W. Veys, \emph{Stringy invariants of normal
surfaces}, J. Alg. Geom. \textbf{13} (2004), 115-141.
\end{thebibliography}
